# Behavior Analysis and Design of Concrete-Filled Steel Circular-Tube Short Columns Subjected to Axial Compression

**Duc-Duy Pham, Phu-Cuong Nguyen***

Faculty of Civil Engineering, Ho Chi Minh City Open University,

Ho Chi Minh City, Vietnam

* cuong.pn@ou.edu.vn; henycuong@gmail.com

**Abstract**

In this paper, a new finite element (FE) model using ABAQUS software was developed to investigate the compressive behavior of Concrete-Filled Steel Circular-Tube (CFSCT) columns. Experimental studies indicated that the confinement offered by the circular steel tube in a CFSCT column increased both the strength and ductility of the filled concrete. Base on the database of 663 test results CFSCT columns under axial compression are collected from the available literature, a formula to determine the lateral confining pressures on concrete. Concrete-Damaged Plasticity Model (CDPM) and parameters are available in ABAQUS are used in the analysis. From results analysis, a proposed formula for predicting ultimate load by determining intensification and diminution for concrete and steel. The proposed formula is then compared with the FE model, the previous study, and the design code current in strength prediction of CFSCT columns under compression. The comparative result shows that the FE model, the proposed formula is more stable and accurate than the previous study and current standards when using material normal or high strength.

## 1. Introduction

A Concrete-Filled Steel Circular-Tube (CFSCT) column is composed of a steel tube filled with concrete, which utilizes the advantages of both concrete and steel tubes. The high compressive strength and stiffness of concrete core in combination with high tensile strength and ductility of steel tube result in many beneficial characteristics for composite section. Therefore, CFSCT columns are considerably more using than RC columns in the high-rise building because of the economy it brings. Moreover, the reason for using the CFST columns is because it substantially increases the workspace area and considerably decreases the cross-section size. Several researchers studied behavior of CFST columns aims to develop the new formulation of CFST columns for design specifications. Although in recent years, many experiments using high-strength materials for both concrete and steel on circular short CFST columns but the comparison of simulation results with experiment is limited. Therefore, the numerical investigation to behavior short circular CFST columns is needed.

In recent years, many researchers have evaluated the characteristics of the CFST column [1-11]. A wide-ranging experimental analysis on composite columns fabricated from normal-high-strength materials for CFST columns by Sakino et al [2], Liew and Xiong [12] and Xiong et al. [13]. With the study of Sakino et al [2], the compressive strength of the concrete ranged from 41MPa to 80MPa whilst the yield stress of steel tubes ranged from 262 MPa to 835 MPa for 114 specimens. Based on the experimental results, the author proposed the calculation to predict the ultimate strength of the composite column. Furthermore, Liew and Xiong [12, 13] investigated the axial performance of concentrically compressed ultra-high-strength CFST short columns. The use of experimental models manufactured with compressive strength of concrete up to 195 MPa and yield stress of steel up to 830 MPa. In their experiments, considering the confining effect of concrete will improve the ultimate strength and when using a small D/t ratio corresponding to ductility was increased. Tang et al. [14], Sakino et al. [2], Gupta et al. [9], Hatzigeorgiou [15], Han et al. [16] and Hu et al. [17] were experimental and analysis, which their purpose is to investigate the development of a more accurate model and the compression behavior of the CFST columns. They focus on the proposed model of concrete confinement with the analysis of the confinement effect in the concrete core of the circular CFST columns. Nevertheless, the hoop stress of the steel tube is the cause of the confinement effects and the different proposed models result in considerably different confinement effects. Therefore, the use of software such as ABAQUS or ANSYS provides a solution for simulating CFST columns through which parameters and effects can be examined.

Addition, the design approaches adopted in paper EC4 [18], AISC [19], ACI [20], CISC [21], DBJ [22] are reviewed and applied to calculate the ultimate strength of the tests columns. Subsequently, the predicted values are compared with the experimental results obtained from the experiments. In general, all the codes somewhat overestimated the capacities except CISC, DBJ present best prediction which $N_{test}/N_u$ difference of about under 5%. EC4, AISC, ACI predicted higher capacity than the experimental results about 6% to 33%; whilst proposed model and design fomula predicted highest about 1.2% than the test.

The main purpose of the paper is to develop a new stress-strain curve for the FE model to predict accurately the static behavior of CFSCT columns by using concrete-damaged plasticity model in ABAQUS. Furthermore, the new models will determine the parameters in the concrete-damaged plasticity model and the stress-train model will be developed for confined concrete. The data collected by the author is used to calibrate the proposed model.

## 2.Finite element modeling of circular CFST column

### *2.1 General description*

The development needs an accurate FE model for predicting the strength and behavior of circular short CFST columns under axial compression by using the finite element ABAQUS. Fig. 1a shows that the symmetric property of the specimens, so during the simulation, it is only necessary to show one-eighth. Since initial imperfections are taken into account in the present model, an eigenvalue buckling analysis is first carried out to provide the lowest buckling mode to be used as the shape of the initial imperfection in subsequent load-deflection nonlinear analysis. The ultimate strength and post-buckling behavior of CFST columns are then obtained from the load-deflection nonlinear analysis.

### *2.2 Element type symmetry and mesh size*

The experimental observations show that the concrete core and steel tube in circular CFST columns don't rotate axial longitudinal. Therefore, Fig.1a shows that CAX4R element (4-node bilinear, reduced integration with hourglass control) is used in the modelling axisymmetric for the steel and concrete components; while Fig.1b is used full 3D model with C3D8R element for both concrete and steel components Previous studies on mesh convergence have been done on circular CFST columns; in which selecting element sizes in lateral and longitudinal directions were chosen D/10 shown in Fig. 1 with D was the diameter of the CFSCT column. Although have the same element size but the number elements of full 3D model are many times higher than axisymmetric model. Moreover, both models provide high accuracy but the analysis time of full 3D model is twenty more than that of axisymmetric model shown in Fig. 2.

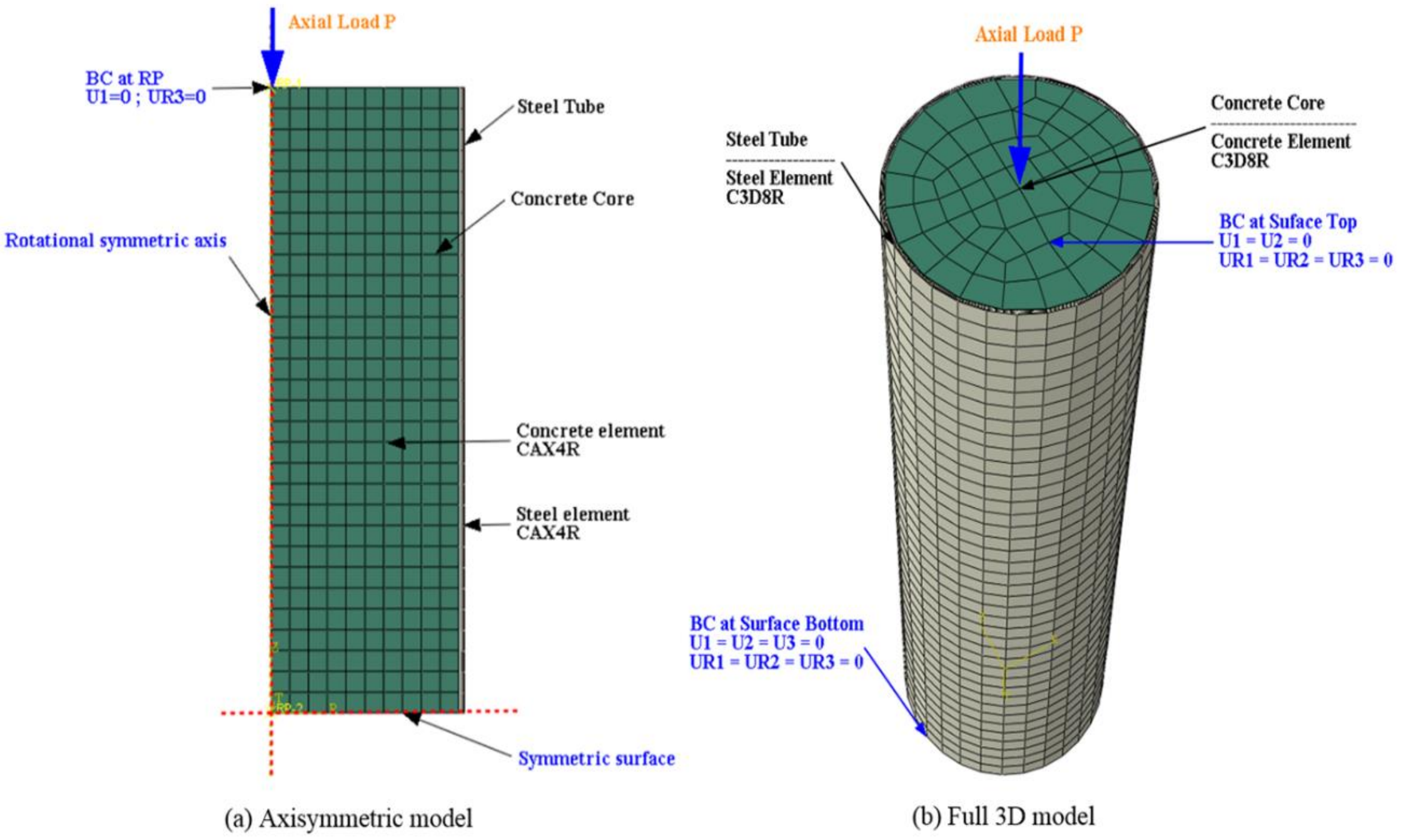

Fig. 1 Modelling of CFST columns

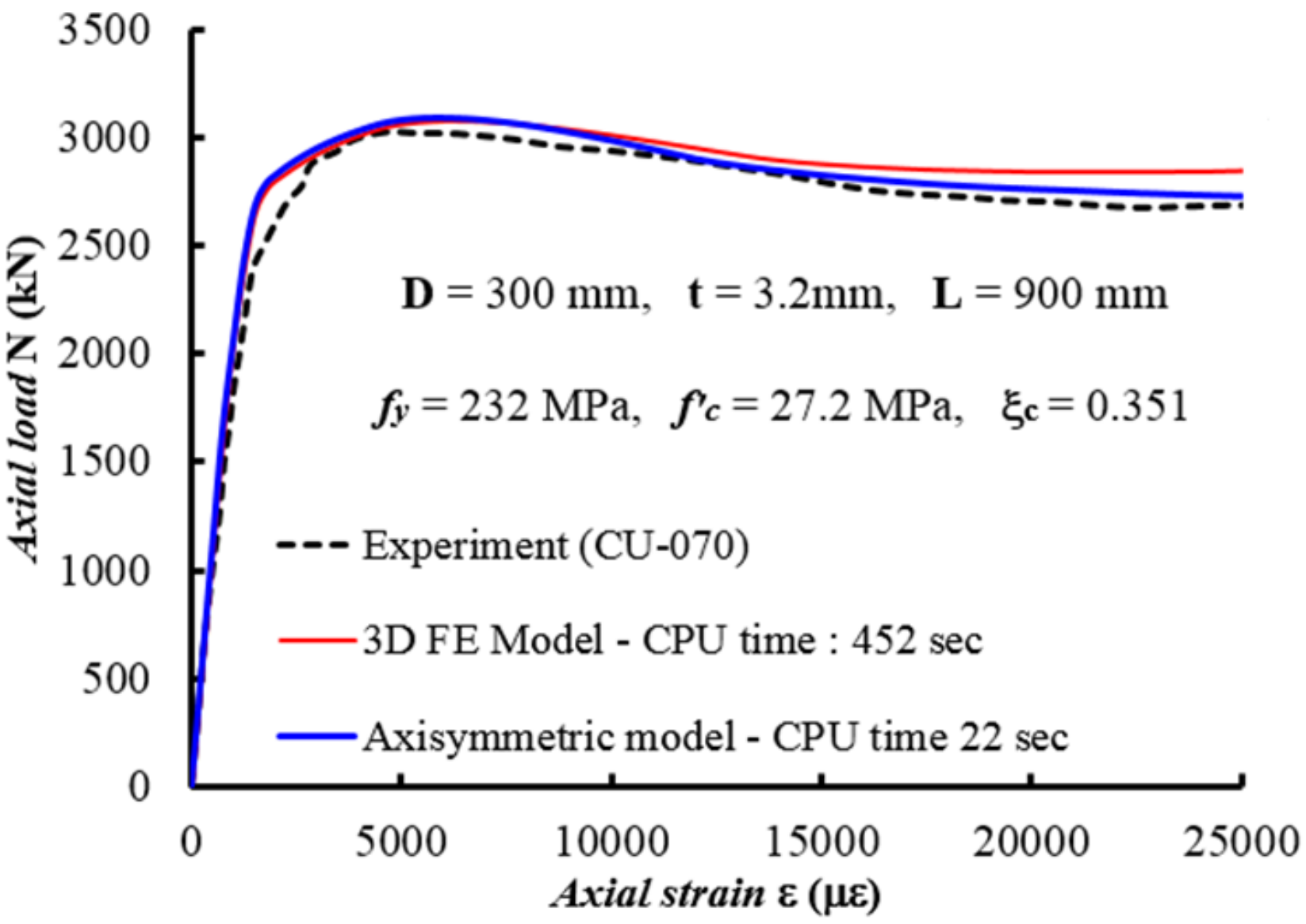

Fig. 2 Comparison between 3D model and the axisymmetric model

*2.3 Contact interaction, loading, boundary*

The surface-to-surface contact in ABAQUS represents the contact interaction of the surface between the steel box and concrete core with *Hard contact* conditions for the *NORMAL BEHAVIOUR option and *Coulomb* friction model for the *TANGENT BEHAVIOUR option. Principally, the coefficient of friction between steel and concrete ranges from 0 to 0.6 in the cases of greasing interface and ungreasing interface, respectively. It was found that there was little effect on the axial resistance when different friction factors were used, but using a smaller friction factor induced a convergent problem with large deformation. Lam and Dai [20] shows that a friction factor of 0.2 or 0.3 is suggested to achieve a quick convergence and to obtain an accurate result. In this study, a friction coefficient of the steel tube and concrete was taken equal to 0.6 as Han et al. [16].

Reference points located at the centroid of the section using a rigid body constraint that is connected to the top and bottom surfaces of CFST columns. With this constraint, the end sections are planar during the analysis and hence there is no need to include the endplates or stiffeners. The loading and boundary conditions are then directly applied to reference points. Due to the symmetry of the specimen, the symmetric boundary conditions also apply to the bottom surfaces at the symmetric planes of the analyze as shown in Fig. 1a. In contrast, the top end is free to displace in the loading direction and the displacement loading is applied at the top reference point.

*2.4 Modeling of initial imperfections*

The initial imperfection of the specimens can be included in the load-deflection analysis using the *IMPERFECTION option available in ABAQUS. In this study, the shape of the local initial imperfections is assumed as the first buckling mode shape obtained from the eigenvalue buckling analysis. It should be noted that the eigenvalue buckling analysis of the composite model is impossible to be used. Thus, the eigenvalue analysis of the corresponding hollow model will be used to obtain the first buckling mode shape as shown in Fig. 3. So the magnitude of the initial imperfection was taken as L/1000 as recommended by Chen et al [23].

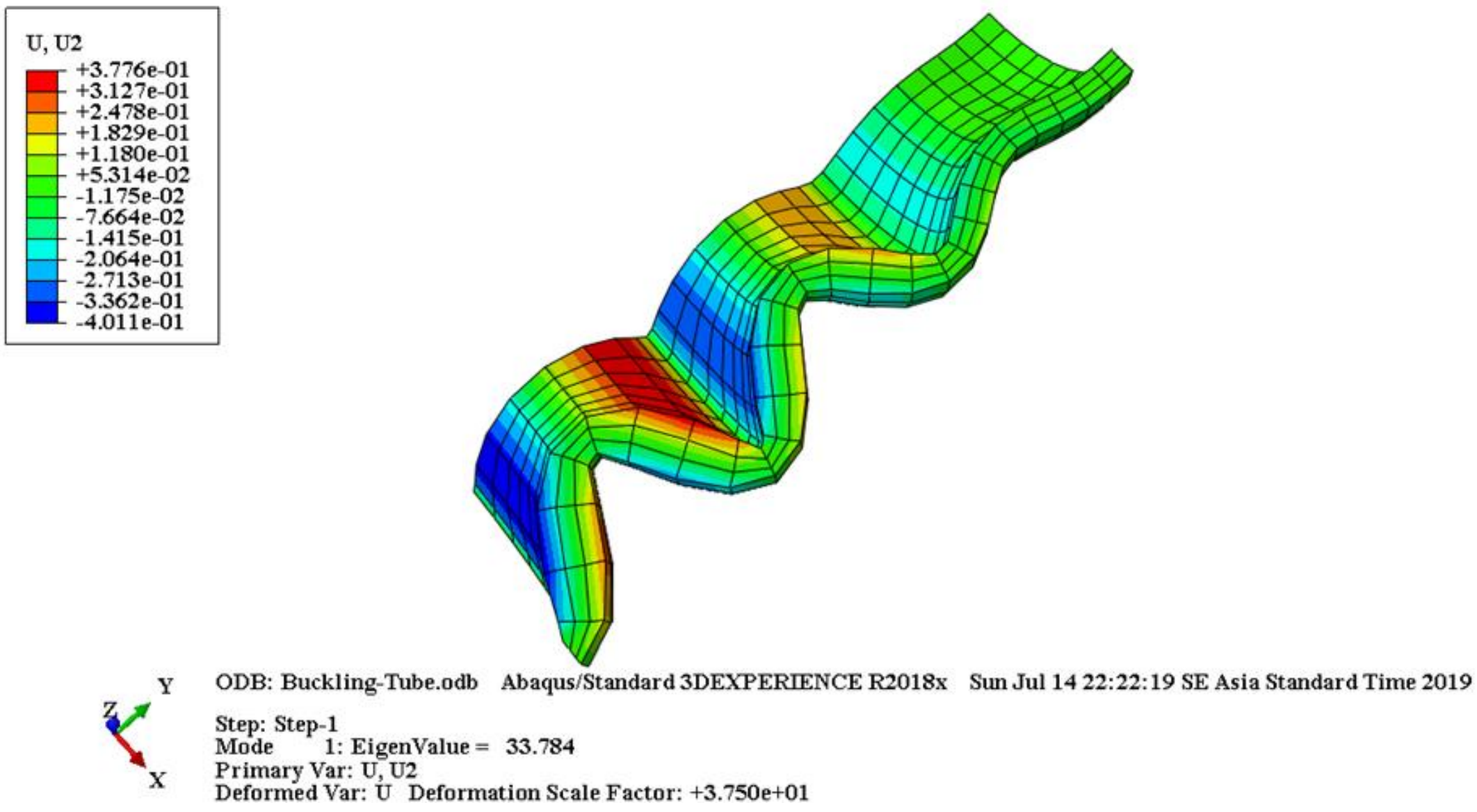

Fig. 3 First buckling mode shape

## 3.Constitutive Models

### *3.1 Constitutive modeling for steel*

Different stress-strain models used for steel tube simulation in CFST column researches such as elastic-perfectly plastic model [17, 24], and elastic-plastic model with linear hardening [25, 26] or multi-linear hardening [16]. The stress-strain relationship model σ - ε proposed by Tao et al. [27] is used in this study as shown in Fig. 4. The model σ - ε of steel proposed with a validity range $f_y$ from 200 MPa to 800 MPa, which is expressed follows:

$$\sigma = \begin{cases} E_s \times \varepsilon & \text{(a)} \\ f_y & \text{(b)} \\ f_u - \left(f_u - f_y\right) \times \left(\dfrac{\varepsilon_u - \varepsilon}{\varepsilon_u - \varepsilon_p}\right)^p & \text{(c)} \\ f_u & \text{(d)} \end{cases} \qquad (1)$$

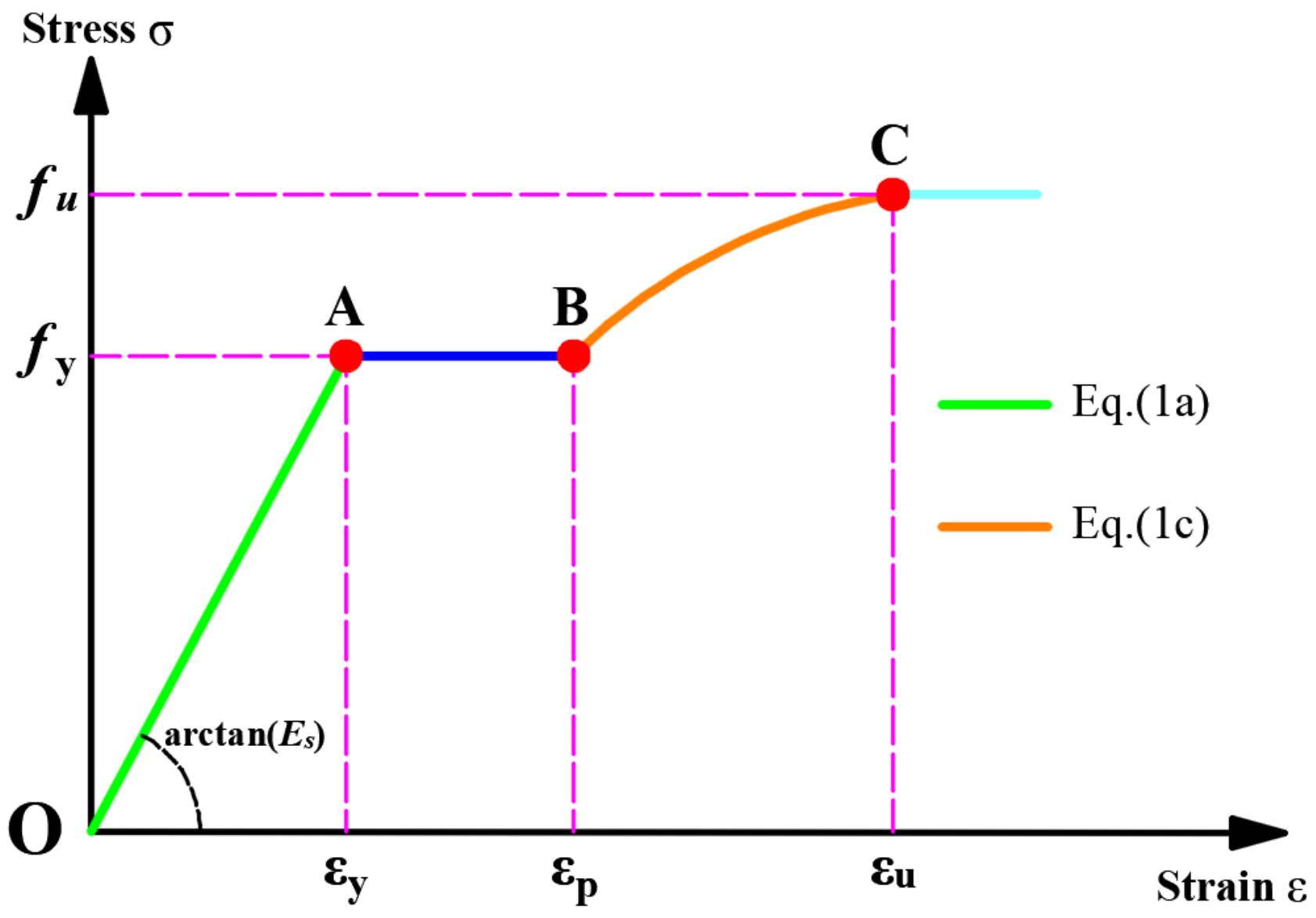

Fig. 4 σ–ε model proposed by Tao et al for structural steel

In which $f_u$ is the ultimate strength; $\varepsilon_y$ is the yield strain, $\varepsilon_y = f_y/E_s$; $\varepsilon_p$ is the strain at the onset of strain hardening. $\varepsilon_u$ is the ultimate strain at the ultimate strength. $p$ is the strain-hardening exponent, which can be determined as follow:

$$p = E_p \times \left(\frac{\varepsilon_u - \varepsilon_p}{f_u - f_u}\right) \qquad (2)$$

in which $E_p$ is the initial modulus of elasticity at the onset of strain hardening and can be taken as $0.02E_s$. $\varepsilon_p$ and $\varepsilon_u$ can be determined as follow:

$$\varepsilon_p = \begin{cases} 15 \times \varepsilon_y & f_y \le 300 \text{ MPa} \\ \left[15 - 0.018 \times \left(f_y - 300\right)\right] \times \varepsilon_y & 300 < f_y \le 800 \text{ MPa} \end{cases} \qquad (3)$$

$$\varepsilon_u = \begin{cases} 100 \times \varepsilon_y & f_y \le 300 \text{ MPa} \\ \left[100 - 0.15 \times \left(f_y - 300\right)\right] \times \varepsilon_y & 300 < f_y \le 800 \text{ MPa} \end{cases} \quad (4)$$

Also, the three parameters are shown in Fig. 4 are yield strength ($f_y$), ultimate strength ($f_u$), and modulus of elasticity ($E_s$), which are required to determine the full-range stress-strain curve.

### *3.2 Concrete damaged plasticity model*

The concrete-damaged plasticity model is used to describe confined concrete for CFST columns and which available in ABAQUS. In this model, material parameters can be determined through formulas. Include the ratio of the second stress invariant on the tensile meridian to that on the compressive meridian ($K_c$), dilation angle ($\psi$), strain hardening – softening rule, the modulus of elasticity ($E_c$), ratio of the compressive strength under biaxial loading to uniaxial compressive strength ($f_{bo}/f_c'$), and fracture energy ($G_f$). Other parameters include flow potential eccentricity (e) and viscosity parameter used default values 0.1 and 0 in ABAQUS respectively.

#### *3.2.1* Determining $K_c$ - $\psi$ - $f_{bo}/f_c'$ - $G_F$

Based on the test data, Papanikolaou and Kappos [28] proposed the following equation to predict the ratio of $f_{bo}/f_c'$

$$\frac{f_{b0}}{f_c'} = \frac{1.5}{\left[f_c'\right]^{0.075}} \quad (5)$$

The ratio of the second stress invariant on the tensile meridian to that on the compressive meridian ($K_c$) is one of the parameters for determining the yield surface of concrete.

According to Yu et al. [29] indicate that the range of $K_c$ varies from 0.5 to 1 and the following equation can be calculated $K_c$:

$$K_c = \frac{5.5\times f_{b0}}{3\times f_c^{'} + 5\times f_{b0}} = 5.5\times\frac{1}{5+2\times\left(f_c^{'}\right)^{0.075}} \tag{6}$$

The dilation angle ($\psi$) is one of the parameters required for ABAQUS to define the yielding potential. In ABAQUS, the dilation angle $\psi$ ranges from $0^0$ to $56^0$. The following equation is proposed by Tao et al. [30] based on regression analysis to determine $\psi$ for circular CFST columns:

$$\psi = \begin{cases} 56.3\times(1-\xi_c) & \xi_c \le 0.5 \\ 6.672\times e^{\frac{7.4}{4.64+\xi_c}} & \xi_c > 0.5 \end{cases} \tag{7}$$

Where the confinement factor $\xi_c$ is expressed as

$$\xi_c = \frac{A_s\times f_y}{A_c\times f_c} \tag{8}$$

The effects of ligament length, rate of loading, and concrete composition on the specific fracture energy $G_f$ and the strain-softening diagram are investigated. The following equation is used to define $G_f$ for tensile concrete by Becq-Giraudonb et al [31]

$$G_f = \left(0.00469\times d_{\max}^2 - 0.5\times d_{\max} + 26\right)\times\left(\frac{f_c}{10}\right)^{0.7}\times 10^{-3}\ \text{(N/mm)} \tag{9}$$

### 3.2.2 Constitutive modeling for concrete

In the initial stages of loading before the growth, the Poisson's ratio of steel tube is greater than that of concrete core. The different lateral expansions of steel tube and concrete core cause a tendency of separation of two materials. Consequently, two materials are assumed to exhibit uniaxial stress states, thus carrying load independently of each other. When the applied load increases and the compressed concrete starts to plasticize. In this situation, the lateral expansion of the concrete reaches maximal value, thus mobilizing the steel tube and efficiently confining the concrete core. For further increase in the load, the steel tube restrains the concrete core and the hoop stresses in the steel become tensile. At this stage and later, the concrete core is subjected to a triaxial stress state induced by radial confining pressure and longitudinal stress, while the steel tube develops a biaxial stress state including longitudinal stress and transverse hoop stress.

The constitutive concrete model given by [30] considering the confinement effects is employed in this axisymmetric analysis. A new proposed material law for the concrete components illustrated in Fig. 5 included three-stages: ascending branch (from Point O to Point A), plateau branch (from Point A-B), descending branch (from Point B to Point C). The material law of the confining concrete suggested by Samani and Attart [32] is considered to compute the nonlinear curve O-A of the stress-strain response, as given by:

$$\frac{\sigma_{Stage\,OA}}{f_c^{'}} = \frac{A\times\left(\frac{\varepsilon}{\varepsilon_{c0}}\right)+B\times\left(\frac{\varepsilon}{\varepsilon_{c0}}\right)^2}{1+\left(A-2\right)\times\left(\frac{\varepsilon}{\varepsilon_{c0}}\right)+\left(B+1\right)\times\left(\frac{\varepsilon}{\varepsilon_{c0}}\right)^2} \quad 0<\varepsilon\leq\varepsilon_{c0} \tag{10}$$

Where $A=\frac{E_c\times\varepsilon_{c0}}{f_c^{'}}$ ; $B=\frac{\left(A-1\right)^2}{0.55}-1$;

The strain at the peak stress of the unconfined concrete is calculated in Eq.(11) proposed by Tasdemir et al. [33] based on the regression analysis of 228 uniaxial compression test specimens with uniaxial compressive strength values ranging from 6 to 105 MPa.

$$\varepsilon_{c0}=\left(-0.067\times\left(f_c^{'}\right)^2+29.9\times f_c^{'}+1053\right)\times10^{-6} \tag{11}$$

After the OA stage, the peak strength value remains constant range $\varepsilon_{c0}$ to $\varepsilon_{cc}$. This stage represents an increment of peak strain and strength of concrete increases due to confinement by the interaction between the steel tube and concrete. The strain at the point

B ($\varepsilon_{cc}$) for the concrete model is determined by the following equation proposed by Xiao et al. [34], which are found to provide an accurate prediction for both normal and high-strength concrete:

$$\frac{\varepsilon_{cc}}{\varepsilon_{c0}} = 1 + 17.4 \times \left( \frac{f_r}{f_c^{'}} \right)^{1.06} \tag{12}$$

In the elastic stage, there is no confining stress ($f_r$). Just before and after the yielding of the steel, the confining stress increases very fast. When once the ultimate strength is reached, confining stress keeps stable or increases very slowly. The procedure for $f_r$ determination is shown in Figure 6. The result of that process is shown in Fig. 7, it shows that the relationship between $f_r$, D/t, $f_y$. It is shown that $f_r$ increases with increasing $f_y$ or decreasing D/t ratio. However, $f_c$' and $\xi_c$ have no significant influence on fr. Based on regression analysis, Eq.(13) is proposed to determine fr for circular columns.

$$f_r = \left[ \frac{1 + 0.03224 \times f_y}{1 + 1.52 \times 10^{-6} \left( f_c^{'} \right)^{-4.5}} \right] \times \exp\left( -0.0212 \times \frac{D}{t} \right) \tag{13}$$

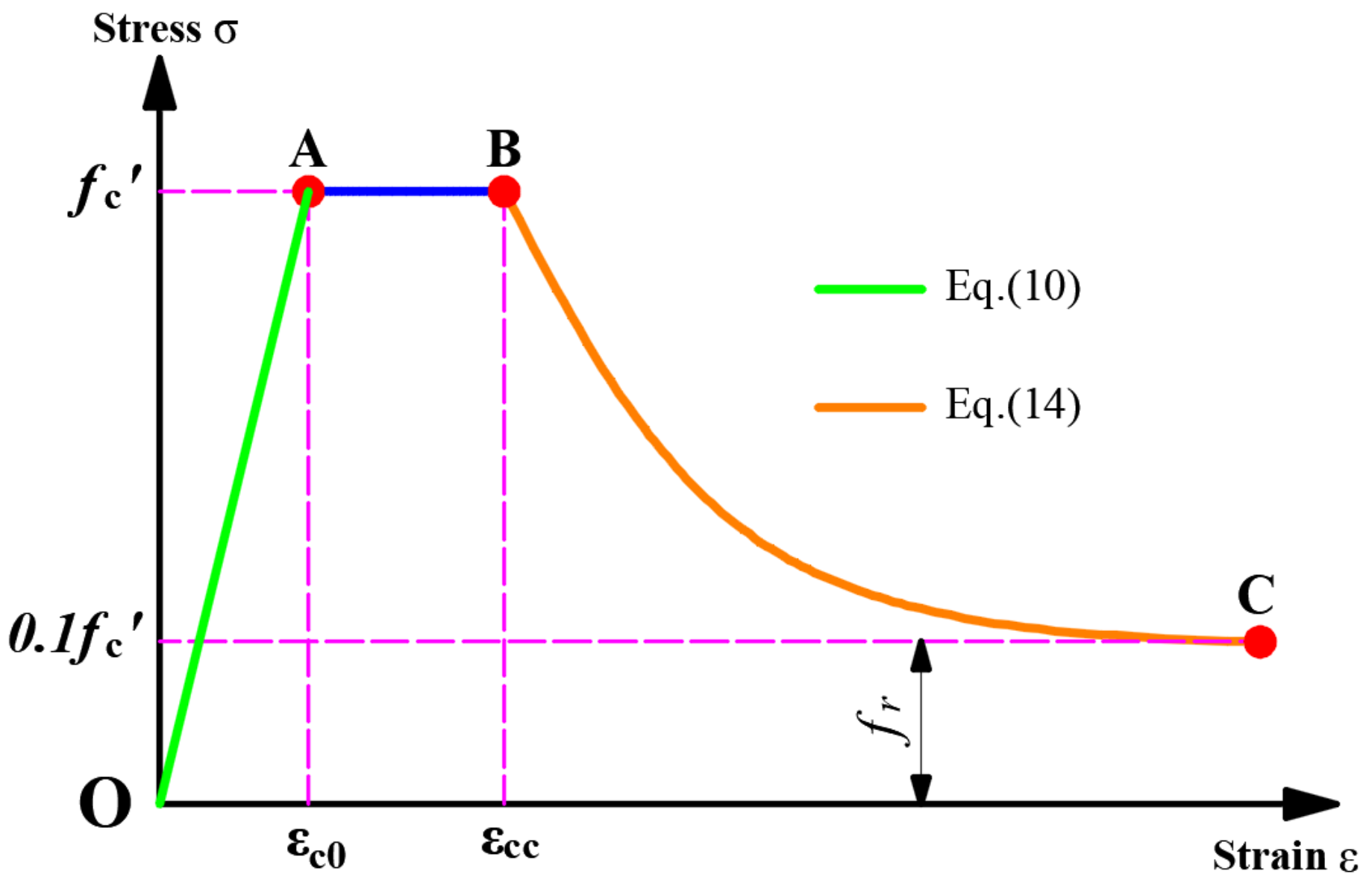

Fig. 5 Stress-strain model for confined concrete in circular CFST columns

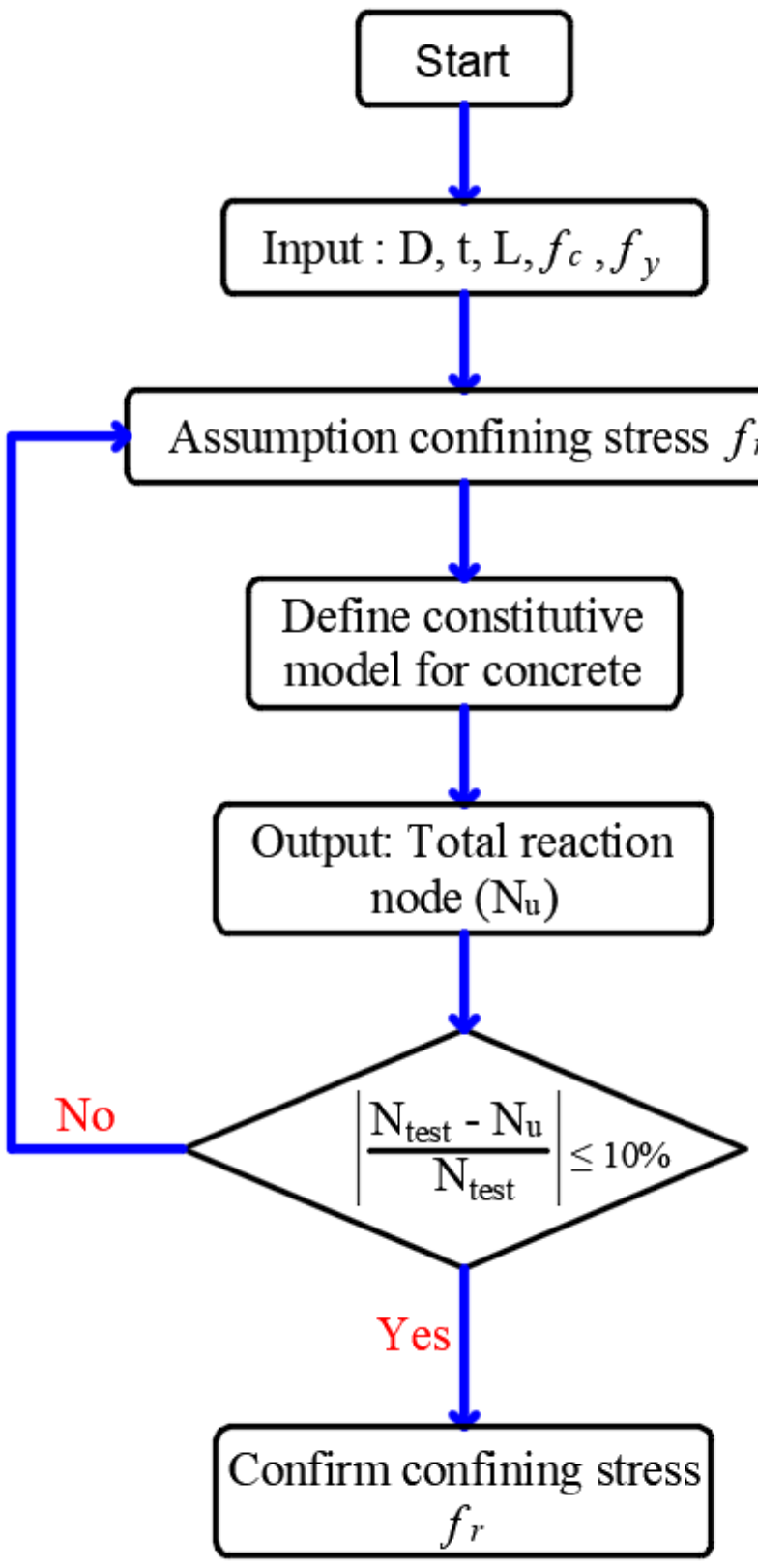

Fig. 6 Flowchart to determine to confine stress ($f_r$) using finite element analysis

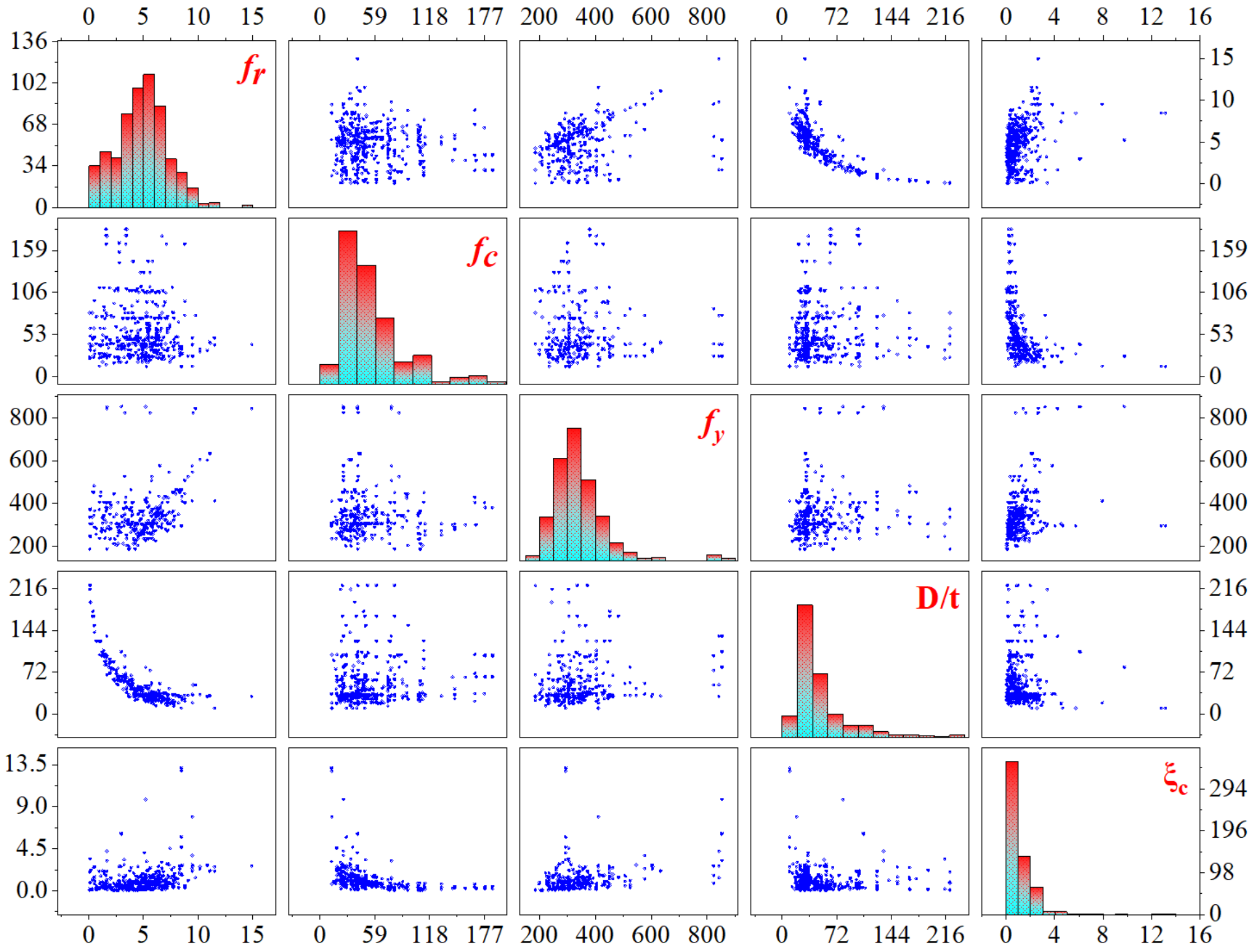

Fig. 7 Evaluate factors affecting $f_r$

An exponential function proposed by Binici [35] is used to represent the descending branch (BC):

$$\sigma_{Stage\ BC} = f_{re} + \left( f_c^{'} - f_{re} \right) \times \exp\left[ -\left( \frac{\varepsilon - \varepsilon_{cc}}{\alpha} \right)^{\beta} \right] \qquad \varepsilon \geq \varepsilon_{cc} \qquad (14)$$

Parameters in the function include the residual stress ($f_r$), the shape of the softening branch ($\alpha$ and $\beta$). These three parameters are proposed by Tao et al. [30].

$$f_{re} = 0.7 \times \left( 1 - e^{-1.38\xi_c} \right) \times f_c^{'} \leq 0.25 \times f_c^{'} \quad (15)$$

$$\alpha = 0.04 - \frac{0.036}{1 + e^{6.08\xi_c - 3.49}}\ ;\ \beta = 1.2 \quad (16)$$

## 4.Verification of the proposed FE Model

Through an extensive literature search, N–ε or N–Δ curves of 663 circulars are collected and used to verify the proposed FE model. Data for the circular columns shown in Table 1 the ranger of different parameters are $f_y$ = 185.7–853 MPa; $f_c'$ = 12.5 – 185.6 MPa; D = 48 – 1020 mm, L/D = 0.8 – 5 and D/t = 10.1 – 220.9. As regards the classification of concrete compressive strength, the following limitations were suggested by Liew and Xiong [36, 37] with Normal strength concrete (NSC) $f_c' \leq 60$ MPa, High strength concrete (HSC) $60 < f_c' \leq 120$ and Ultra-high strength concrete (UHSC) $f_c' \geq 120$ MPa. Database descriptions collected are summarized through the histograms, which are illustrated in Fig. 8. It can be seen that a large number of tests using NSC ($f_c' \leq 60$ MPa about 61.24%) and normal strength steel ($f_y \leq 460$ MPa about 91.70%) have been conducted on the CFST columns. Only a small number of tests have been carried out on high,ultra-high-strength steel ($f_y > 460$ MPa about 8.3%) , HSC ($60 < f_c' \leq 120$ MPa about 33.63%) and UHSC ($f_c' \geq 120$ MPa about 5.13%). As observed, almost tested specimens have been fabricated with large steel tubes thickness (D/t ≤ 60 about 77.38%), it less used much with a small thickness (D/t > 60 about 22.62%).

The compressive strength of concrete obtained from the compressive tests on standard cylinders of 150mm × 300mm ($f_{cyl,150}$) is used for Eq (5-14). Experimental data collected in table 1 using the compressive strength of concrete $f_c$ were varied such as 150mm×300mm ($f_{cyl,150}$), 100mm×300mm ($f_{cyl,100}$), 100mm×100mm ($f_{cu,100}$), 150mm×150mm($f_{cyl,150}$). Consequently, the values of the compressive strength $f_c$ in the collected experiments were converted into $f_{cyl,150}$ by using the conversion factors proposed by Skazlic et al. [38] for UHSC and by Yi et al. [39] for NSC and HSC as follows:

$$f_c' = \begin{cases} f_{cyl,150} = \dfrac{f_{cyl,100}}{1.03} = 0.88\times f_{cu,150} = 0.82\times f_{cu,100} & \text{with NSC} \\ f_{cyl,150} = \dfrac{f_{cyl,100}}{1.04} = 0.98\times f_{cu,150} = 0.92\times f_{cu,100} & \text{with HSC} \\ f_{cyl,150} = 0.95 f_{cyl,100} & \text{with UHSC} \end{cases} \quad (17)$$

Table 1 Experimental data collected

| **Source** | ***N*** | ***D*** (mm) | ***t*** (mm) | ***D/t*** | ***L/D*** | $f_y$ (MPa) | $f_c'$ (MPa) |
|---|---|---|---|---|---|---|---|
| [40] | 05 | 140.10 – 168.30 | 4.47 – 9.68 | 14.47 – 37.65 | 2.89 – 4.83 | 270 – 302 | 23.40 – 33.20 |
| [41] | 14 | 76.43 – 152.63 | 1.68 – 4.93 | 29.5 – 48.45 | 2.0 | 363.3 – 633.42 | 43.44 – 20.9 |
| [42] | 12 | 168.30 – 169.30 | 2.60 – 5.00 | 33.80 – 64.62 | 1.81 | 200.2 – 338.10 | 18.20 – 37.10 |
| [43] | 03 | 218.30 – 219.25 | 6.10 – 6.50 | 33.60 – 36.24 | 4.3 | 302.00 | 37.07** |
| [44] | 10 | 108.00 – 108.60 | 4.50 – 4.60 | 23.60 – 24.00 | 2.20 – 4.12 | 271.9 – 409.60 | 28.60 – 30.70 |
| [45] | 30 | 76.50 – 300.00 | 1.50 – 4.50 | 24.00 – 100.00 | 2.00 – 5.00 | 232.30 – 602.70 | 9.90 – 48.35 |
| [46] | 17 | 166.00 – 320.00 | 5.00 – 7.00 | 33.20 – 45.71 | 0.8 – 3.98 | 249.90 – 274.40 | 26.60 – 46.57 |
| [47] | 19 | 121.00 – 320.00 | 5.00 – 12.00 | 10.10 – 45.71 | 0.8 – 3.98 | 250.10 – 294.20 | 9.20 – 52.96** |
| [48] | 12 | 96.00 – 273.00 | 2.00 – 12.00 | 10.10 – 102.000 | 1.70 – 4.22 | 235.20 – 410.60 | 15.70 – 40.28 |
| [47] | 21 | 96.00 – 273.00 | 2.00 – 12.00 | 10.10 – 102.00 | 4.00 – 4.69 | 235.40 – 410.90 | 11.90 – 46.88** |
| [49] | 05 | 130.60 – 134.30 | 2.40 – 6.20 | 21.60 – 43.04 | 2.00 – 2.02 | 235 | 17.40 – 26.60 |
| [50] | 06 | 100.00 – 101.80 | 0.50 – 5.70 | 17.90 – 192.31 | 2.00 | 244.00 – 320.00 | 18.00 – 37.40 |
| [51] | 06 | 108.00 | 4.00 | 27.00 – 34.60 | 3.00 – 3.14 | 324.00 – 339.11 | 29.00 – 34.03** |
| [50] | 10 | 100.00 – 101.80 | 0.50 – 5.70 | 17.86 – 192. | 311.96 – 2.00 | 244.19 – 319.70 | 17.95 – 37.36 |
| [48] | 05 | 48.00 – 165.00 | 0.70 – 4.50 | 13.71 – 214.29 | 3.20 – 4.00 | 223.00 – 304.00 | 22.56 – 33.44 |
| [52] | 04 | 86.49 – 89.27 | 2.74 – 4.05 | 22.02 – 32.66 | 3.02 – 3.12 | 226.70 | 30.20 – 48.00 |
| [53] | 02 | 114.30 | 3.50 | 25.40 | 2.00 | 339.00 – 350.00 | 33.40 |
| [53] | 12 | 174.00 – 179.00 | 3.00 – 9.00 | 19.78 – 58.00 | 2.01 – 2.07 | 248.50 – 283.31 | 22.16 – 45.70 |
| [54] | 10 | 159.00 – 1020.00 | 5.07 – 13.25 | 32.36 – 105.81 | 3.00 | 291.40 – 381.50 | 15.00 – 46.00 |
| [55] | 09 | 76.00 – 101.70 | 2.20 – 2.40 | 34.55 – | 3.5 | 380.00 – 390.00 | 57 |
| [56] | 01 | 323.90 | 5.60 | 89.41 | 3.09 | 443.90 | 92.30 |
| [57] | 02 | 152 | 1.70 | 89.41 | 3.29 | 270.00 | 73.00 |
| [58] | 01 | 165.20 | 4.17 | 39.62 | 4.00 | 358.70 | 40.90 |
| [59] | 01 | 165.20 | 4.50 | 36.71 | 4.00 | 413.84 | 40.89 |
| [60] | 15 | 165.00 – 190.00 | 0.86 – 2.82 | 58.51 – 220.93 | 3.47 – 3.52 | 185.70 – 363.30 | 41.00 – 108.00 |
| [61] | 09 | 111.30 – 165.70 | 2.00 – 6.30 | 23.65 – 56.78 | 2.94 – 3.07 | 309.50 – 482.50 | 60.75 |
| [62] | 05 | 165.00 – 190.00 | 0.86 – 2.82 | 58.51 – 220.93 | 3.47 – 3.52 | 185.70 – 363.30 | 41.00 – 48.30 |
| [63] | 06 | 190.00 | 1.11 | 171.17 | 3.45 – 3.49 | 203.10 | 94.70 – 110.30 |
| [24] | 02 | 140.80 – 141.40 | 3.00 – 6.50 | 21.75 – 46.93 | 4.49 – 4.51 | 285.00 – 313.00 | 23.81 – 28.18 |
| [57, 64] | 31 | 101.30 – 265.00 | 1.50 – 5.40 | 23.65 – 85.50 | 2.87 – 3.93 | 264.90 – 357.70 | 18.40 – 42.70 |
| [65] | 16 | 108.00 – 133.00 | 1.00 – 7.00 | 18.14 – 125.00 | 3.5 | 232.00 – 429.00 | 77.40 – 106.02 |
| [66] | 26 | 101.30 – 318.50 | 3.03 – 10.37 | 30.69 – 33.63 | 2.99 – 3.01 | 335.00 – 452.00 | 23.20 – 52.20 |
| [1] | 15 | 165.00 – 190.00 | 0.86 – 2.82 | 58.51 – 220.93 | 3.47 – 3.52 | 185.70 – 363.30 | 41.00 108.00 |
| [67] | 03 | 157.50 – 158.70 | 0.90 – 2.14 | 73.69 – 176.33 | 2.84 – 2.86 | 221.00 – 308.00 | 18.70 |
| [3] | 13 | 114.09 – 115.04 | 3.75 – 5.02 | 22.91 – 30.53 | 2.61 – 2.63 | 343.00 – 365.00 | 31.40 – 104.90** |
| [2] | 36 | 108.00 – 450.00 | 2.96 – 6.47 | 16.69 – 152.03 | 3.00 | 279.00 – 853.00 | 25.40 – 85.10 |
| [68] | 26 | 60.00 – 250.00 | 1.87 – 2.00 | 30.00 – 133.69 | 3.00 | 282.00 – 404.00 | 85.20 – 90.00 |
| [69] | 17 | 165.00 – 219.00 | 2.43 – 4.78 | 45.85 – 67.90 | 2.97 – 3.09 | 350.00 | 34.10 – 68.00 |
| [70] | 04 | 100.00 | 1.90 | 52.63 | 3.00 | 404.00 | 121.60 |
| [71] | 09 | 114.30 | 3.35 – 6.00 | 19.05 – 34.12 | 3.00 | 287.33 – 342.95 | 32.68 – 105.45 |
| [72] | 09 | 150.00 – 450.00 | 3.20 | 46.88 – 140.63 | 4.11 – 4.22 | 265.00 | 25.40 |
| [73] | 03 | 111.64 – 113.64 | 1.90 – 3.64 | 31.22 – 58.76 | 3.52 – 3.58 | 259.60 – 261.30 | 47.80 – 56.70 |
| [74] | 07 | 75.84 – 76.21 | 2.48 – 3.31 | 23.02 – 30.69 | 3.94 – 3.96 | 278.00 – 305.00 | 145.00 |
| [75] | 36 | 129.00 – 133.00 | 3.00 – 5.00 | 26.60 – 43.00 | 3.00 | 306.00 | 53.70 – 76.40 |
| [76] | 06 | 114.30 | 2.74 – 5.90 | 19.37 – 41.72 | 2.62 | 235.00 – 355.00 | 56.20 – 107.20 |
| [77] | 01 | 165 | 2.37 | 69.62 | 3.73 | 287.50 | 30.20** |
| [78] | 03 | 557.65 – 559.40 | 16.52 – 16.54 | 33.76 – 33.82 | 1.78 – 1.79 | 546 | 31.70** |
| [79] | 12 | 114.30 – 219.10 | 3.60 – 10.00 | 18.14 – 43.82 | 2.19 – 2.74 | 380.00 – 428.00 | 51.60 – 193.30 |
| [13] | 16 | 114.30 – 219.10 | 3.60 – 10.00 | 21.91 – 43.82 | 2.19 – 2.74 | 300.00 – 428.00 | 51.60 – 193.30 |
| [80] | 35 | 153.00 – 477.00 | 1.54 – 11.36 | 41.52 – 102.20 | 2.00 | 290.00 – 345.00 | 73.20 |
| [81] | 09 | 127.00 – 203.00 | 4.20 – 9.60 | 21.15 – 30.24 | 3.61 – 3.94 | 349.00 – 427.00 | 84.20 |
| [82] | 09 | 107.90 – 114.90 | 2.05 – 8.03 | 55.46 – 13.44 | 2.98 – 3.17 | 251.80 – 304.30 | 59.00 – 130.80 |
| [83] | 10 | 100.00 – 168.30 | 2.50 – 3.00 | 33.33 – 60.11 | 1.78 – 3.00 | 317.80 – 445.52 | 33.39 – 94.68 |
| [84] | 07 | 200.00 | 6.00 | 33.33 | 3.48 | 451.00 | 38.00 – 112.10 |
| [85] | 18 | 140.00 – 216.30 | 4.50 – 8.20 | 26.38 – 48.07 | 3.00 | 371.90 – 462.90 | 28.00 – 52.00 |
| [47, 51] | 02 | 166.00 | 5.00 | 33.20 | 4.28 | 277.34 | 39.06 – 42.24 |
| [86, 87] | 12 | 100.00 – 200.00 | 3.00 | 33.33 – 66.67 | 3.00 | 303.50 | 46.80 |

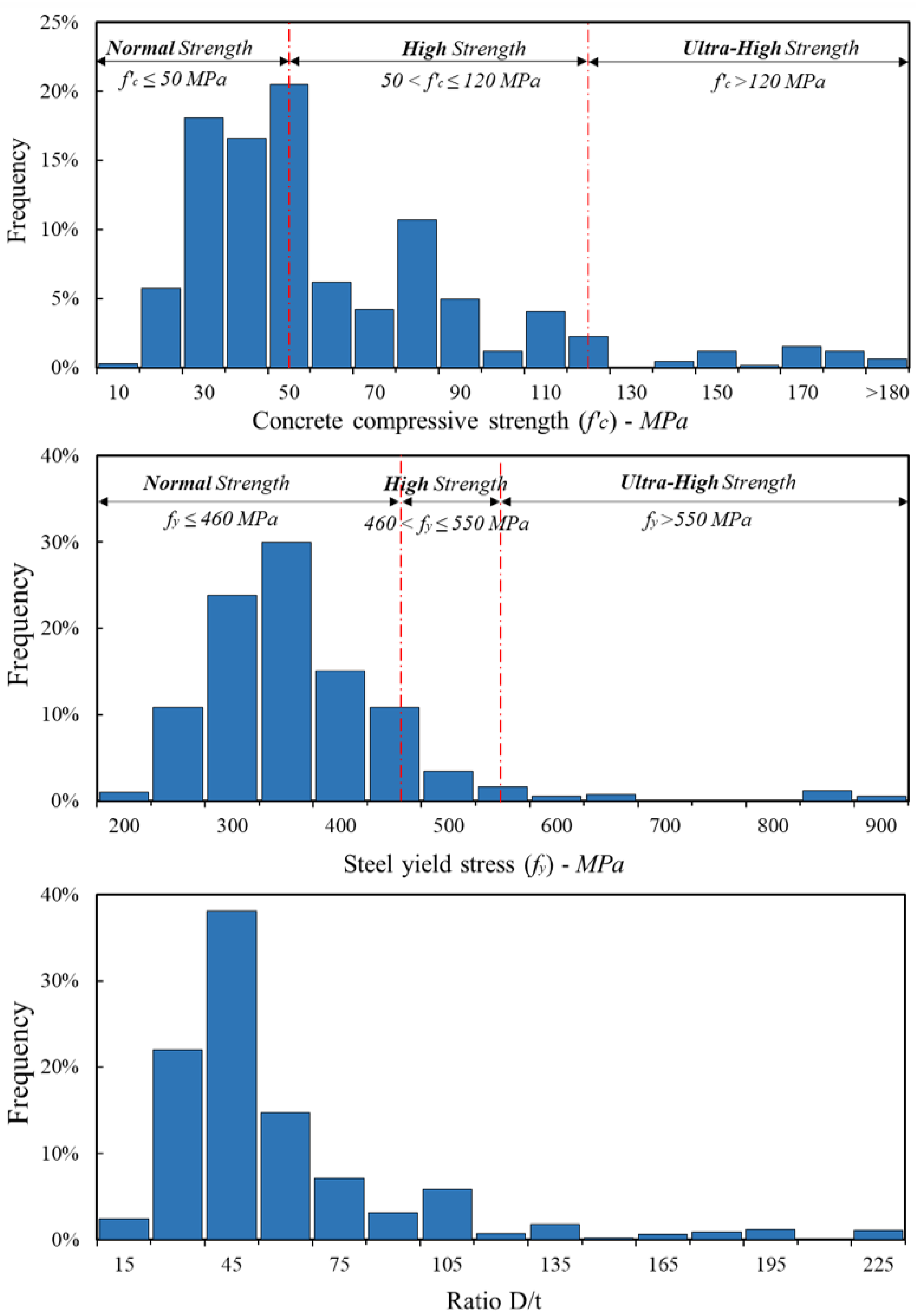

Fig. 8 Distribution of the database

In this study, both FE models are Han [4] and Hu [17] were compared with the proposed FE model. The predicted ultimate strength ($N_{uc}$) from the proposed FE model, Han et al model, and Hu et al are compared with the measured ultimate strength ($N_{exp}$) as shown in Table 2. Table 2 shows both the Mean value (Mean), Standard Deviation (STD), and Coefficient of Variation (CoV) of the ratio of $N_{ue}/N_{uc}$ specimens with different cross-sections.

Table 2 Comparison results of FE predictions with measured ultimate strength.

| FE model | Mean | STD | CoV |
|---|---|---|---|
| Proposed model (Abaqus) | 1.011 | 0.055 | 0.054 |
| Han et al [4] | 0.975 | 0.090 | 0.092 |
| Hu et al [17] | 0.980 | 0.129 | 0.132 |

The mean values of $N_{ue}/N_{uc}$ are 1.011, 0.975 and 0.980; whilst the standard deviations, coefficient of variation 0.055, 0.090, 0.129 and 0.054, 0.092, 0.132 of $N_{ue}/N_{uc}$ for proposed model, Han et al model, Hu et al model. As can be seen, predictions obtained from the proposed FE model with reasonable accuracy, which is superior to Han et al model and Hu et al model in predicting ultimate strength. Fig. 9 show that, with normal strength concrete ($f_c' \leq 60$ MPa), the prediction Han et al, Hsuan The Hu et al are too low and overestimated the ultimate load when compared with the experiment. The predictions of the two models above are reasonable when compared to the experiment in range high ($60 < f_c' \leq 120$ MPa) and ultra-high-strength concrete ($f_c' > 120$ MPa). The prediction Han et al, Hsuan The Hu et al are overestimated safely compared with experiments for both normal and ultra-high-strength steel. The data collection is large, so it can be seen the model proposes a fair evaluation when compared to the experiment for all materials.

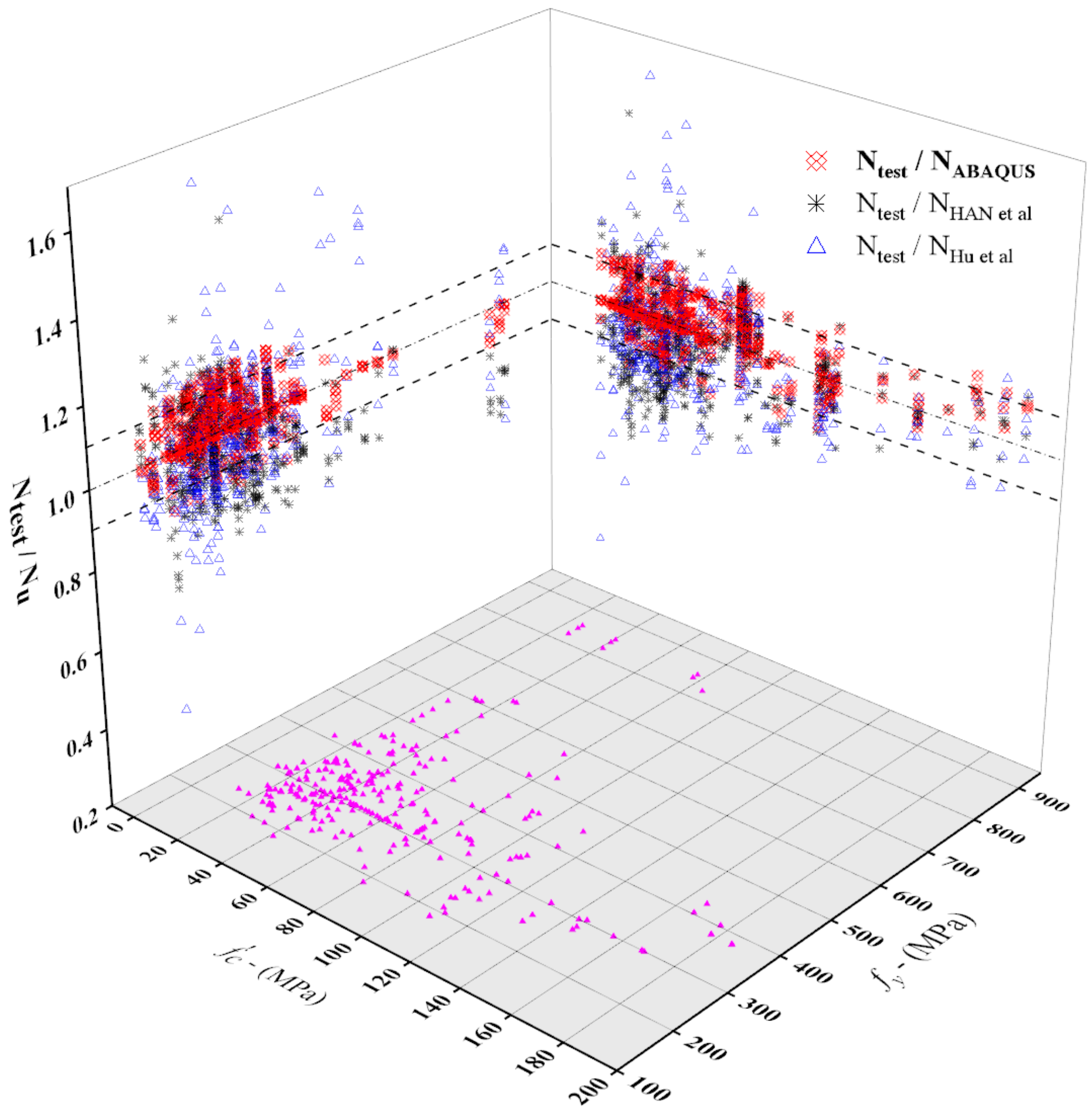

Fig. 9 The effect of strength material on the models

Fig. 10, Fig. 12 shows the prediction is compared with Han et al, Hu et al, and the current FE models for specimens with thin-walled tubes and high-strength concrete. The results show that for circular thin-wall and high strength columns, the proposed model more accurately than the other two models in predicting the strain-softening branch. Two specimens CU-070, CU-100 in Fig. 10 for the predicted and experimental axial load-strain were tested by Hsuan The Hu et al [17]. It can be observed that the current model

closely predicts the initial axial stiffness of the specimen before attaining a loading level of about 80% for the ultimate load. After, local buckling on the test lasts until the load peak is reached, which can be seen in the test. the axil bearing capacity, the columns decrease slightly after the peak load owing to a slight confinement effect applied by the steel tube. Other than Fig. 10 using normal concrete, Fig. 12 using high concrete and thin-walled tubes for specimens CC2-1, CB4-1 was tested by Han et al [68]. The predicted initial axial stiffness of the column is slightly higher than the experimental one. This is likely due to the uncertainty of the actual concrete stiffness and strength as the average concrete compressive strength was used in the simulation. However, when after the post-peak of columns decreases rapidly because the expansion of high strength concrete is very small.

For the two curves shown in Fig. 11, the limitation of the Hu et al model is the ratio D/t < 20 and D/t > 150 so cannot be simulated. This is to make the FE model better than the Hu et al model because there is no limit on the D/t ratio. The figure demonstrates is tested by Sakino et al [2] that these two curves are almost identical up to the ultimate load. With a small D/t = 16.61 for specimens CC8-A-8, softening branches tend to increase steadily because due to effective confinement good for thick tube steel. With a large D/t ratio for specimens CC4-D-4-1 that prediction to increase the ultimate load. After the stage, it can see that the prediction model is not close to the experiment because the local buckling for steel tubes so the stiffness decreases sharply.

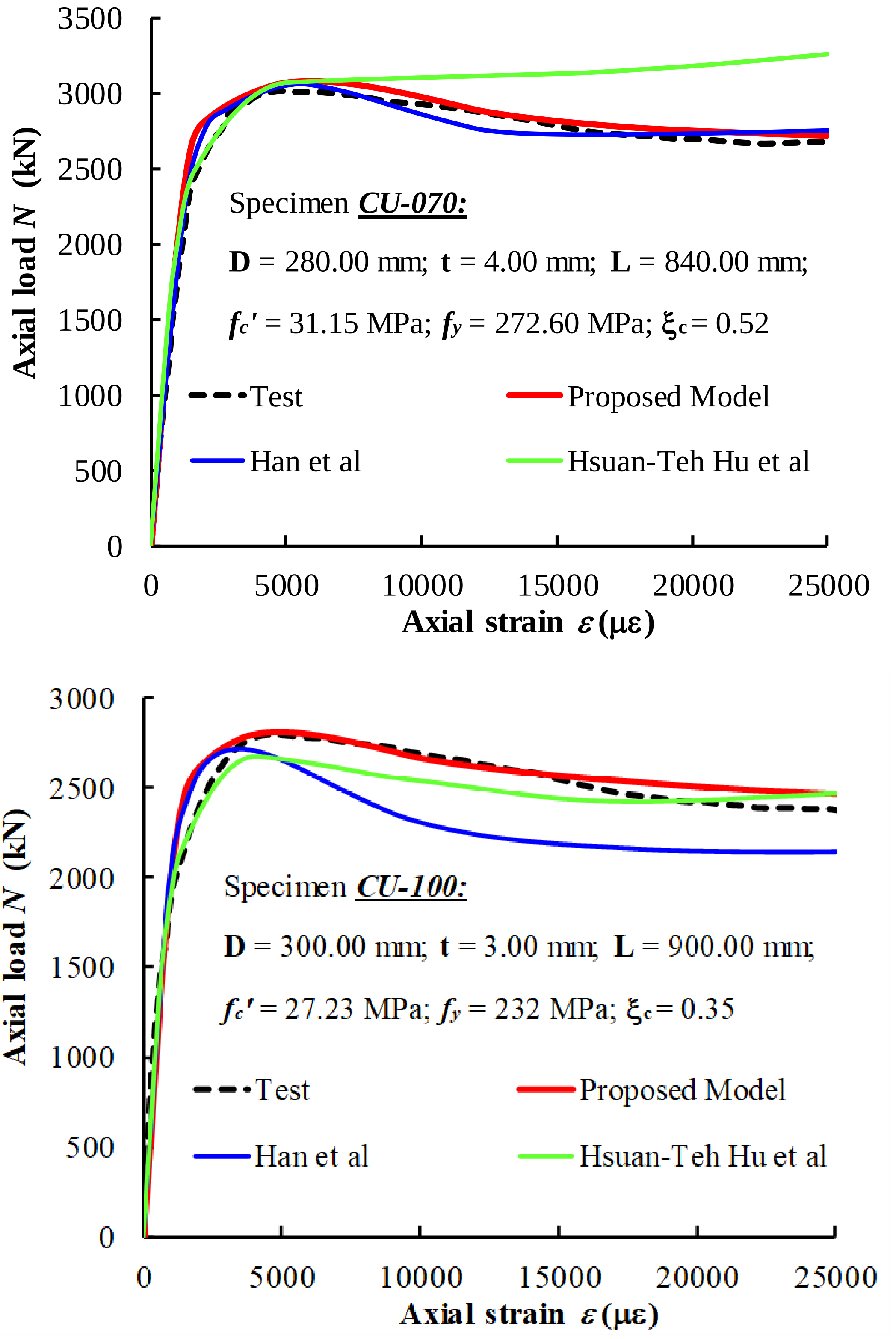

Fig. 10 Comparison between predicted and measured N–ε curves for normal specimens

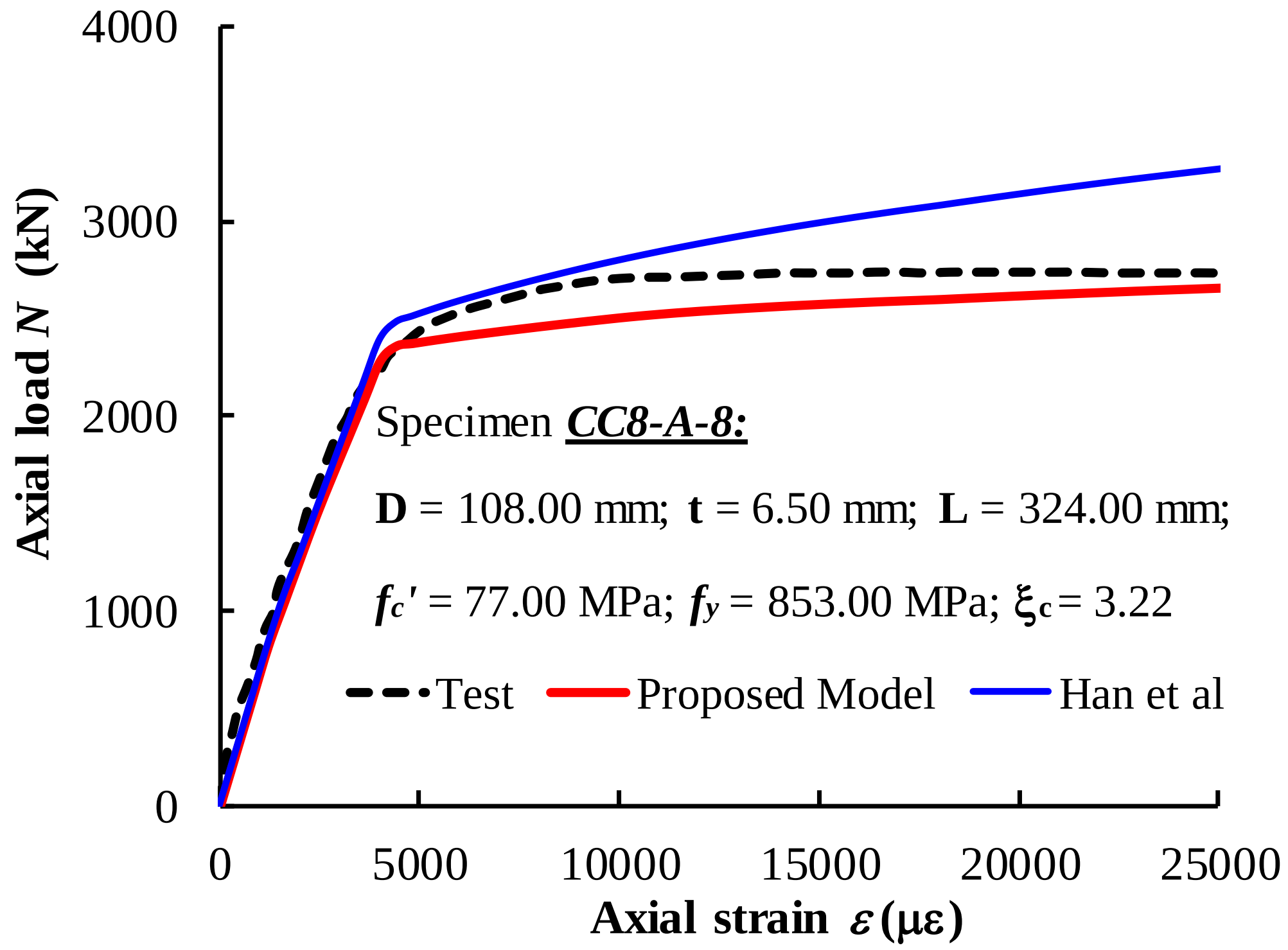

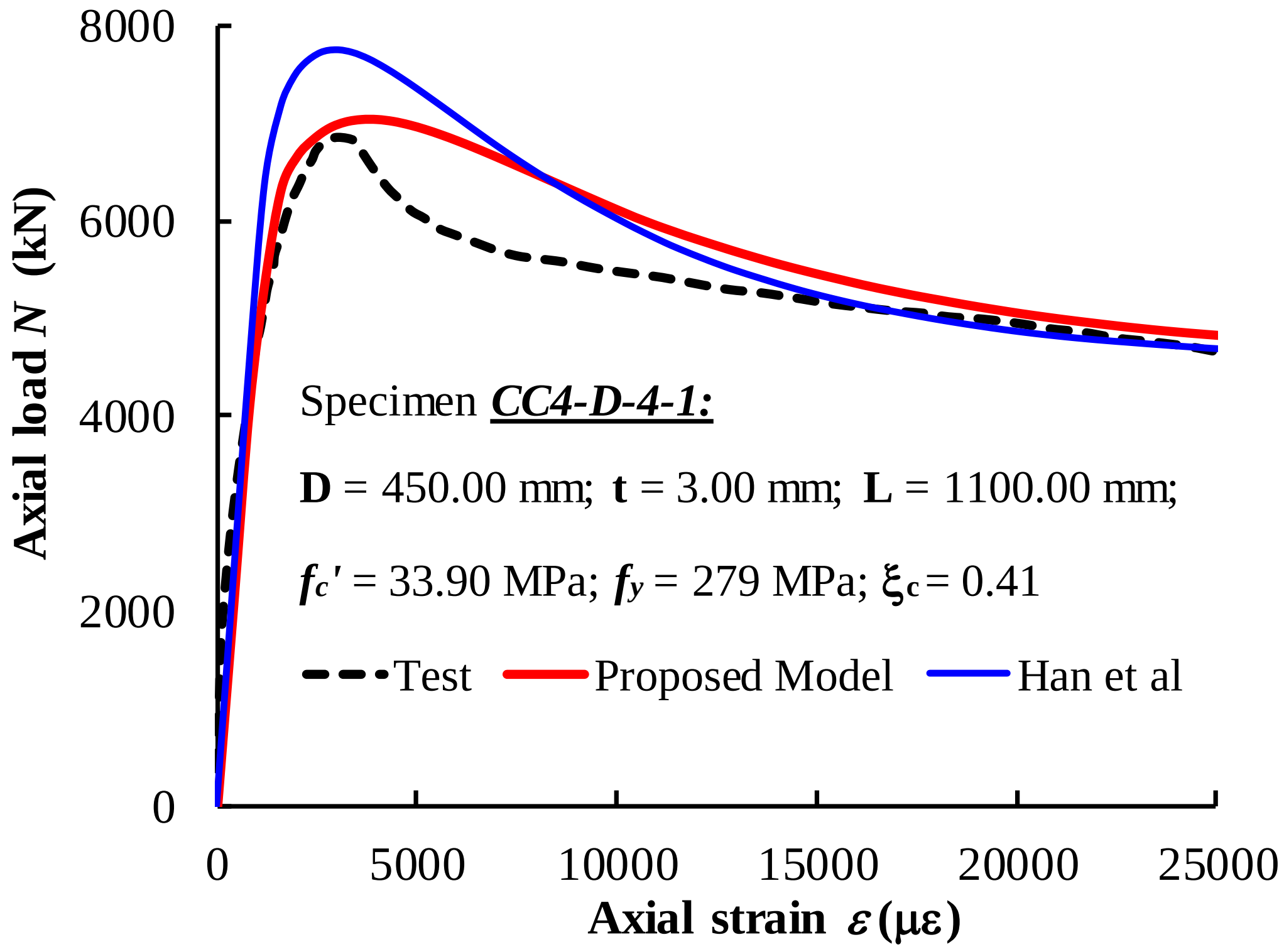

Fig. 11 Comparison between predicted and measured N–ε curves for CFST columns with $21.7 \leq D/t \leq 150$ ratio

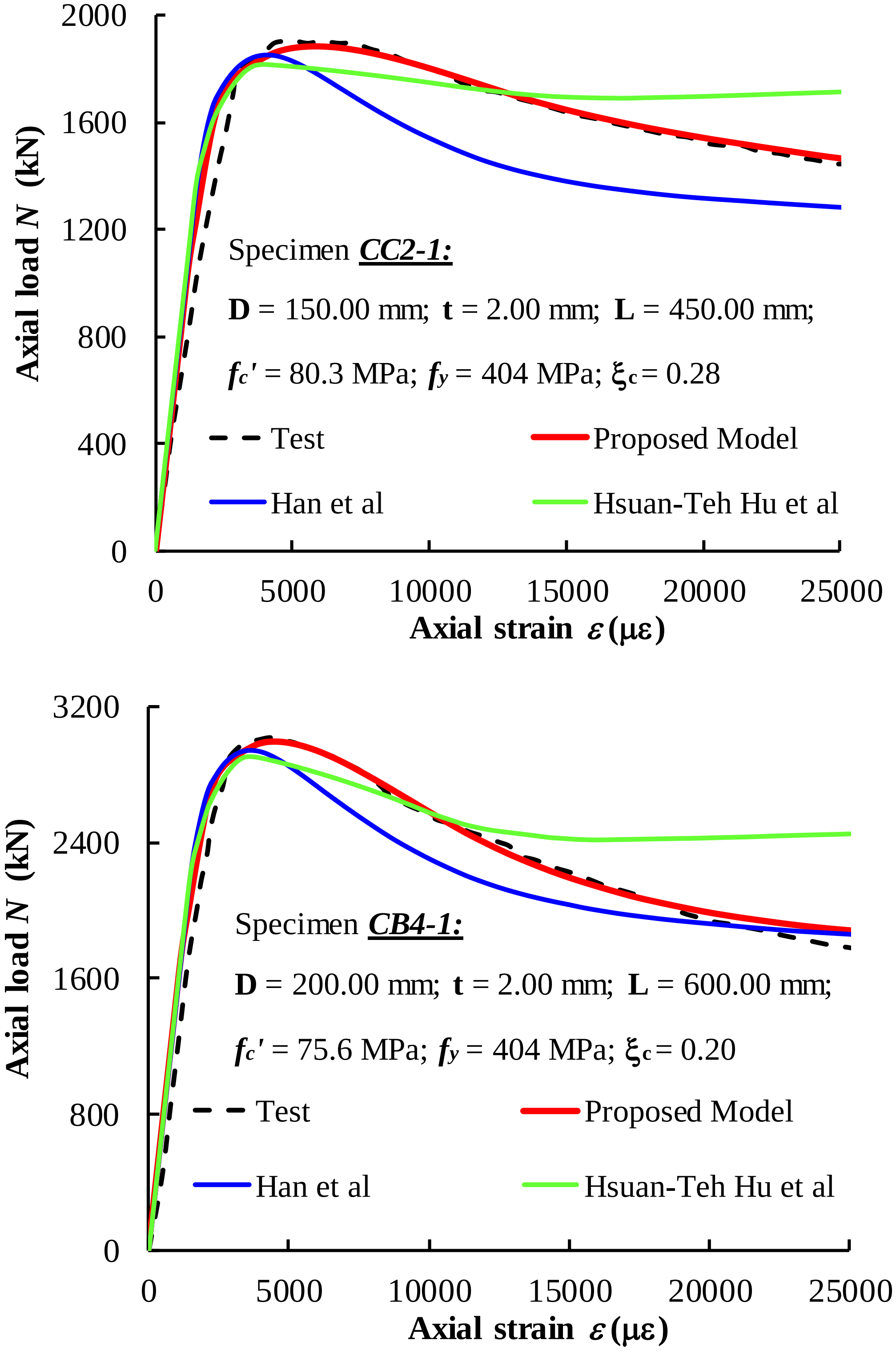

Fig. 12 Comparison between predicted and measured N–ε curves for CFST columns with thin-walled and high strength

## 5. Assessment of current design codes, proposed previous formula and a proposed design for CFST columns

### *5.1 Evaluation of current design codes*

The ultimate loads will be calculated from international standards and these values were consequently compared to the test experimental observation results. The international standards mentioned in the paper include EC4 [18], ANSI/AISC [19], ACI 318R [20], CISC (2007) [21], DBJ 13-51-2010 [22], which respectively condensed to $N_{EC4}$, $N_{AISC}$, $N_{ACI}$ and $N_{CISC}$. The formulas used to predict the ultimate load from five selected codes for the axial compression column CFST are shown in Fig. 13 and Table 4.

As is well know, the application in the use of high strength concrete is not mentioned in these present standards. Consequently, research is needed to apply high-strength materials relevance of current standards for CFST columns. Besides, the prediction of the ultimate load of current design standards is usually based on the combination force-resistance of the core concrete tube steel. With different design codes, the contribution is also different. Table 3 shows the results of comparing $N_{test}/N_u$ ratio between design code and test results. Overall, the values of the Mean of all predictions are higher than the experiment, which highest at 33% (AISC) and lowest at 3% (EC4). In terms of a number test specimen, the limitation specified in the design code, it is still limited to predict the ultimate load on the test specimens. In addition to the two design codes ACI (663/663), CISC (663/663) predict for all test specimen, the remaining standards still cannot be applied as EC4 (472/663), AISC (330/663), DBJ (279/663).

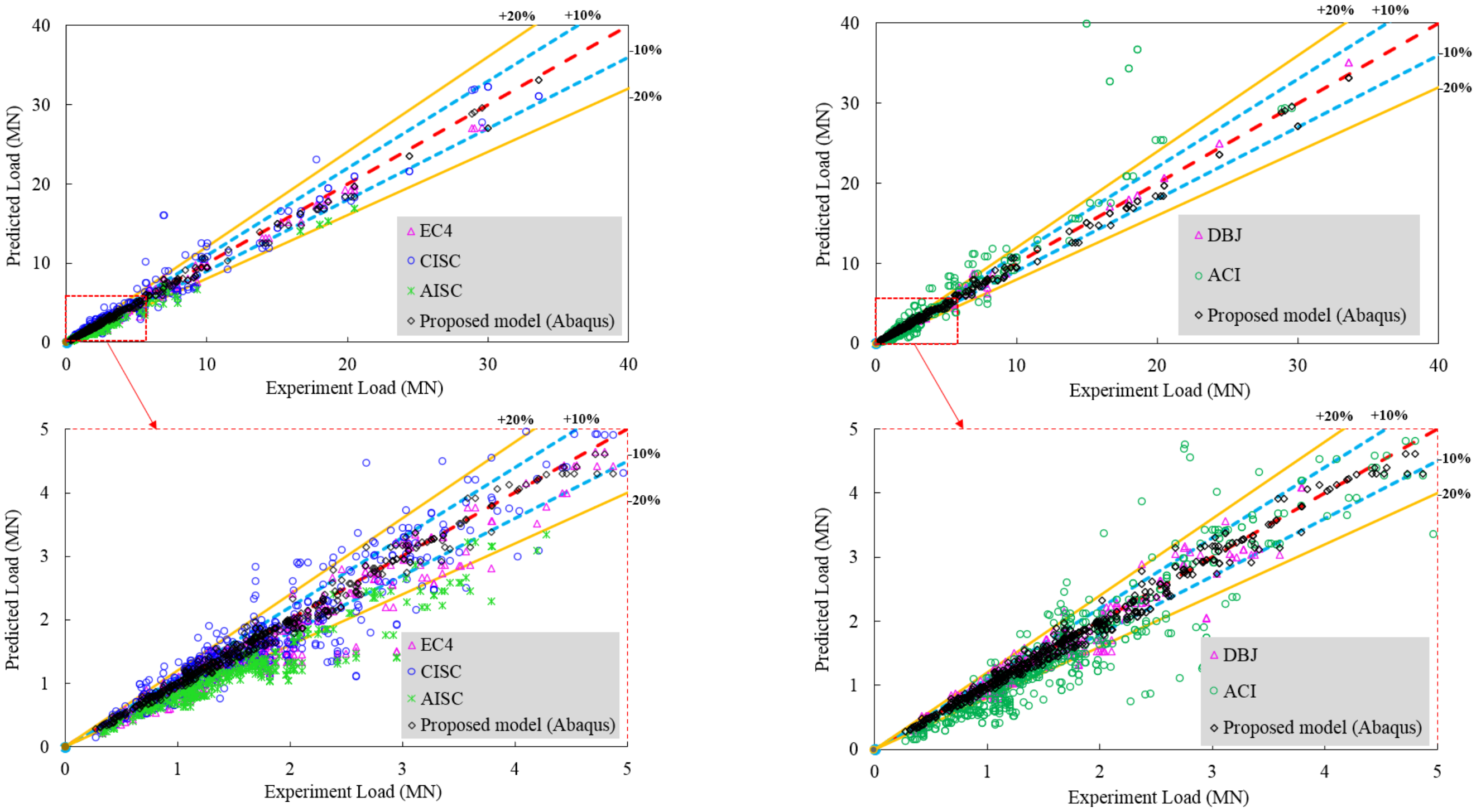

Fig. 13 Comparison between the proposed model and current design codes

Table 3 Evaluation of models in the prediction of current design codes

| FE model | N | Mean | STD | CoV |
|---|---|---|---|---|
| Proposed | 663/663 | 1.011 | 0.055 | 0.054 |
| Eurocode 4 [18] | 472/663 | 1.13 | 0.17 | 0.15 |
| ANSI/AISC [19] | 330/663 | 1.33 | 0.19 | 0.14 |
| ACI 318 [20] | 663/663 | 1.18 | 0.39 | 0.33 |
| CISC [21] | 663/663 | 1.06 | 0.21 | 0.20 |
| DBJ 13-51-2010 [22] | 279/663 | 1.03 | 0.12 | 0.12 |

The prediction $N_{EC4}$ is overestimated the ultimate load. Specifically, the mean value is 13.2% and the values STD and CoV are 0.17 and 0.15 respectively. The results of this prediction show that with $N_{EC4}$ give not more accurate predictions the ultimate load for the confinement effects in circular CFST short columns. AISC shows that with the Mean difference of about 33.45% and CoV and STD are 0.14 and 0.19 respectively when compared to the test results. thereby showing that this method prediction highly overrated the ultimate load. Like above design codes, ACI predicts that ultimate load is considerably different with average value of about 18%. Also, STD and CoV values are the highest among codes, respectively, 0.39,0.33. Likewise, $N_{EC4}$, $N_{AISC}$, and $N_{ACI}$, both $N_{CISC}$, $N_{DBJ}$ underestimated the ultimate loads. But ratio $N_{test}/N_u$ of both codes are the lowest among design codes about 1.06, 1.03.

Table 4 Strength prediction methods and related limitations for design codes

| **Design codes** | **Prediction of ultimate strength** | **Limitations** |
|---|---|---|
| EC4 2004 [18] | $N_{EC4} = \eta_a A_s f_y + A_c f_c \left(1 + \eta_s \frac{t f_y}{D f_{ck}}\right)$ (18) <br> $\eta_c = 4.9 - 18.5\bar{\lambda} + 17\bar{\lambda}^2 \ (\eta_c \geq 0)$; <br> $\eta_a = 0.25(3 + 2\bar{\lambda}) \ (\eta_a \leq 1)$ <br> $\bar{\lambda} = \sqrt{\frac{N_{pl,Rk}}{N_{cr}}}$; $N_{pl,Rk} = f_y A_s + 0.85 f_c A_c$ <br> $N_{cr} = \frac{\pi^2 (EI)_{eff}}{l^2}$; $(EI)_{eff} = E_s I_s + K_e E_{c2} I_c$ | $D/t \leq \sqrt{8E_s / f_y}$ <br> $f_c^{'} \geq 17.2MPa$ |
| AISC 360-10 2010 [19] | $N_{AISC} = \begin{cases} P_{0,AISC}\left[0.658^{\left(\frac{P_{0,AISC}}{P_e}\right)}\right], & (P_e > 0.44P_{0,AISC}) \\ 0.877P_e, & (P_e < 0.44P_{0,AISC}) \end{cases}$ (19) <br> $P_{0,AISC} = 0.95 f_c A_c + f_y A_s$ <br> $(EI)_{eff} = E_s I_s + C_3 E_{c1} I_c$; $C_3 = 0.6 + 2\left(\frac{A_s}{A_s + A_c}\right)$ | $D/t \leq 0.15\frac{E_s}{f_y}$ <br> $f_y \leq 525\ MPa$ <br> $21 \leq f_c^{'} \leq 70\ MPa$ |
| CISC 2007 [21] | $N_{CISC} = (\tau A_s f_y + \tau' 0.85 A_c f_c)(1 + \lambda^{3.6})^{-0.556}$ (20) <br> $\tau = \begin{cases} (1 + \rho + \rho^2)^{-0.5} & L/D < 25 \\ 1.0 & L/D \geq 25 \end{cases}$ <br> $\tau = \begin{cases} 1 + \left(\frac{25\rho^2\tau}{D/t}\right)\left(\frac{f_y}{0.8 f_c}\right) & L/D < 25 \\ 1.0 & L/D < 25 \end{cases}$ <br> $\rho = 0.02(25 - L/D)$ | |

| | | | |
|---|---|---|---|
| | $\lambda = \sqrt{\dfrac{\tau A_s f_y + \tau' 0.85 A_c f_c}{\pi^2 \left( E_s I_s + \dfrac{0.6 E_c I_c}{C_{fs} / C_f} \right) / (KL)^2}}$ | | |
| DBJ 13-51-2010 [22] | $N_{DBJ} = f_{sc}(A_s + A_c)$<br>$f_{sc} = f_{ck}(1.14 + 1.02\xi);\ \xi = \dfrac{f_y A_s}{f_{ck} A_c}$ | (21) | $D/t \le 150 \dfrac{235}{f_y}$<br>$235 \le f_y \le 420\ MPa$<br>$24 \le f_c' \le 70\ MPa$ |
| ACI 318 [20] | $N_{ACI} = A_s f_y + 0.85 A_c f_c$ | (22) | $D/t \le \sqrt{8 \dfrac{E_s}{f_y}}$<br>$f_c' \ge 17.2\ MPa$ |

*5.2 Evaluation of proposed previous formula*

The formulas to predict the ultimate strength ($N_u$) of circular CFST columns include O'Shea and Bridge [1], Yu et al [88], Liu et al [89], Sun [90], Zhong and Miao [44], Guo [91], De Oliveira et al [92]. Moreover, all the above models were directly developed for circular CFST columns. Table 6 shows the sequence in each analytical model to predict $N_u$ and the confined peak stress ($f_{cc}$) of circular CFST columns. The predicted results of the analytical model are compared to the experimental ultimate loads ($N_{u,\ test}$) shows in Fig. 14, Table 5.

Table 5 Evaluation of models in prediction of the proposed previous formula

| FE model | N | Mean | STD | CoV |
|---|---|---|---|---|
| Proposed model (Abaqus) | 663/663 | 1.011 | 0.055 | 0.054 |
| O'Shea and Bridge [1] | 580/663 | 0.77 | 0.22 | 0.28 |
| Yu et al [88] | 141/663 | 1.01 | 0.11 | 0.11 |
| Liu et al [89] | 619/663 | 0.98 | 0.12 | 0.13 |
| Sun [90] | 663/663 | 0.93 | 0.45 | 0.48 |
| Zhong and Miao [44] | 663/663 | 1.12 | 0.21 | 0.18 |
| Guo [91] | 553/663 | 0.81 | 0.10 | 0.12 |
| De Oliveira et al [92] | 629/663 | 1.27 | 0.21 | 0.17 |

Most formulas for predicting the ultimate load of authors apply all specimens collected in Table 1, but the exception of formula Yu applies only very little to specimens collected. Table 5 shows that the models O'Shea and Bridge, Liu et al, Sun, Guo have overestimated the ultimate load compared to the experiment, in which very high Mean values of the models are 0.77, 0.98, 0.93, 0.81 respectively. On the other hand, the model's Yu et al, Zhong and Miao, De Oliveira et al, Wang et al have underestimated the ultimate load compared to the experiment, shown through Mean values respectively, 1.01, 1.12, 1.27, 1.12.

Table 6 The method proposed previously for prediction ultimate load and related limitations

| **Proposed formula design** | **Prediction of ultimate strength** | **Limitations** |
|---|---|---|
| O'Shea and Bridge [1] | $N_u = \sigma_{cp} A_c + A_s f_y$ (23)<br>• $f_c \le 50$ MPa :<br>$\sigma_{cp} = f_c \left( -1.228 + 2.172\sqrt{1 + \frac{7.46 f_l}{f_c}} - 2\frac{p}{f_c} \right)$<br>• $50 \le f_c \le 100$ MPa<br>$\frac{\sigma_{cp}}{f_c} = \left( \frac{p}{f_c} + 1 \right)^k$<br>$p = p_{yield}\left(0.7 - \sqrt{\frac{f_c}{f_y}}\right)\left(\frac{10}{3}\right)$ $p_{yield} = \frac{2 f_y t}{D - 2t}$<br>$k = 1.25\left(1 + 0.062\frac{p}{f_c}\right)(f_c)^{-0.21}$<br>$f_l = 0.558\sqrt{f_c}$ | $D/t \le 200$ |
| Yu et al [88] | $N_u = f_{cc} A_c$ (24)<br>$f_{cc} = (1.14 + 1.34\xi) f_c$ | $235 \le f_y \le 345\ MPa$<br>$30 \le f_c^{'} \le 60\ MPa$<br>$0.2 \le \xi \le 2$ |
| Liu et al [89] | $N_u = \sigma_v A_s + \sigma_{cp} A_c$ (25)<br>$\sigma_v = 0.61 f_y;\ \sigma_h = 0.54 f_y$<br>$\sigma_{cp} = f_c + 4.1\sigma_r;\ \sigma_r = \frac{2t\sigma_h}{D-2t} = \frac{1.08tf_y}{D}$ | |
| Sun [90] | $N_u = f_{cc} A_c$ (26)<br>$f_{cc} = f_c\left(1 + 8.2\frac{(D/t-1)}{(D/t-2)^2}\frac{f_y}{f_c}\right)$ | |
| Zhong and Miao [44] | $N_u = N_c + N_s$ (27)<br>$N_s = \frac{\mu' + 2}{\left[3\left(\mu'^2 + \mu' + 1\right)\right]^{0.5}} f_y A_s$<br>$N_c = \left(f_c + 4p_0\right) A_c$<br>$\mu' = -\frac{1}{2} - \frac{1}{2(\xi+1)};\ \xi = \frac{f_y A_s}{f_c A_c}$ | |
| Guo [91] | $N_u = f_{cc} A_c$ (28)<br>$f_{cc} = f_c\left(1 + \sqrt{\xi} + 1.1\xi\right)$ | $\xi \le 1.7$ |
| De Oliveira et al [92] | $N_u = \left(A_c f_c + A_s f_y\right)\lambda$ (29)<br>$\begin{cases} \lambda = 1 & when\ L/D \le 3 \\ \lambda = -0.18\ln\left(\frac{L}{D}\right) & when\ L/D > 3 \end{cases}$ | $1 \le L/D \le 10$ |

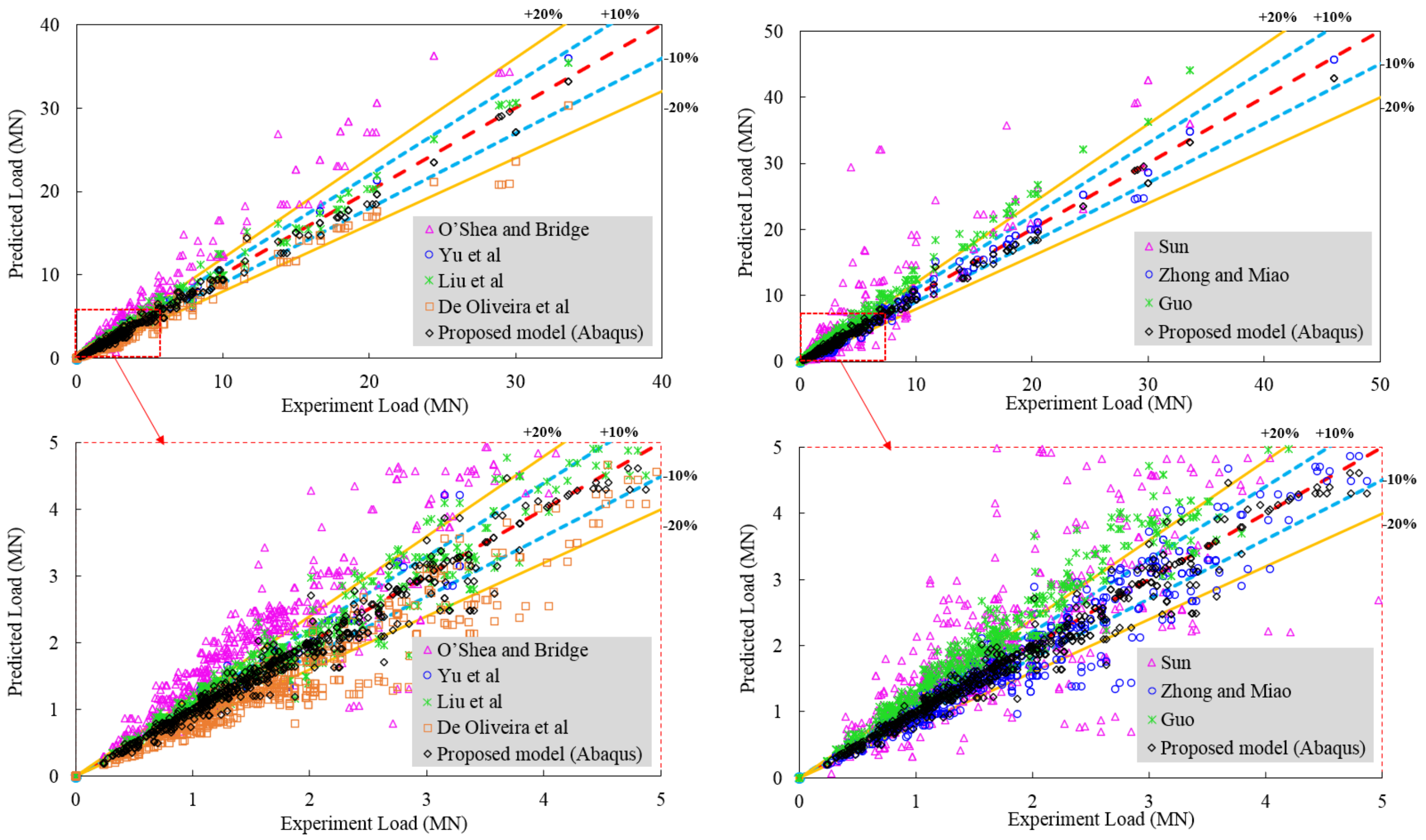

Fig. 14 Comparison between the proposed previous formula and current proposed model

## 5.3 <u>Proposed design</u>

The confined concrete core in triaxial stress and the steel tube is in a biaxial when circular CFST columns under axial compression. The confinement to the concrete core can lead to an increased capacity compressive, however, local buckling of the steel tube can decrease the ultimate load. Through this mechanism can consider contributions for concrete core ($N_c$) and steel tubes ($N_s$) are a combined model. The proposal is as follows to predict the ultimate strength of the circular CFST columns.

$$N_u = N_c + N_s = \eta_c \times \left[ A_c \times f_c^{'} \right] + \eta_s \times \left[ A_s \times f_y \right] \qquad (30)$$

Where $\eta_c$ is an intensification factor and $\eta_c \geq 1$. This factor due to the confinement effect concrete core to increase the compressive strength of concrete. $\eta_s$ is a diminution factor and $\eta_s \leq 1$ due to consider initial local buckling of the steel tubes or residual stress. The cross-section areas of the core concrete ($A_c$) and the steel tubes ($A_s$). It cannot be measured the loads carried by core concrete ($N_c$) and steel tube ($N_s$), so it can be determined from the numerical analysis for the circular CFST columns. Based on Eq. (30), the factors of $\eta_c$ and $\eta_s$ can be determined as:

$$\eta_c = \frac{N_c}{A_c \times f_c^{'}};\ \eta_s = \frac{N_s}{A_s \times f_y} \tag{31}$$

$\eta_s$ is mainly affected by the $\alpha_s = A_s/A_c$ ratio of the cross-section through by D/t ratio and yield stress of the steel tube. In general, $\eta_s$ increases with increasing $f_y$ or decreasing D/t ratio, as depicted in Fig. 15. When $f_y$ increases or D/t ratio decreases, the concrete is under increased confinement. However, the ratio of hoop tensile stress to the yield stress of the steel tube decreases, leading to increased $\eta_s$.

Nonlinear regression is performed and Eq. (32) is proposed to predict $\eta_s$ for circular CFST columns, it shows that stability and high accuracy with standard deviation $R^2$ = 0.8224 shown in Fig. 16.

$$\eta_s = \left[1.923 - 1.229 \times \ln\left(0.003 f_y\right)\right] \times \left(\alpha_s \times \frac{f_c^{'}}{f_y}\right)^{0.47} \tag{32}$$

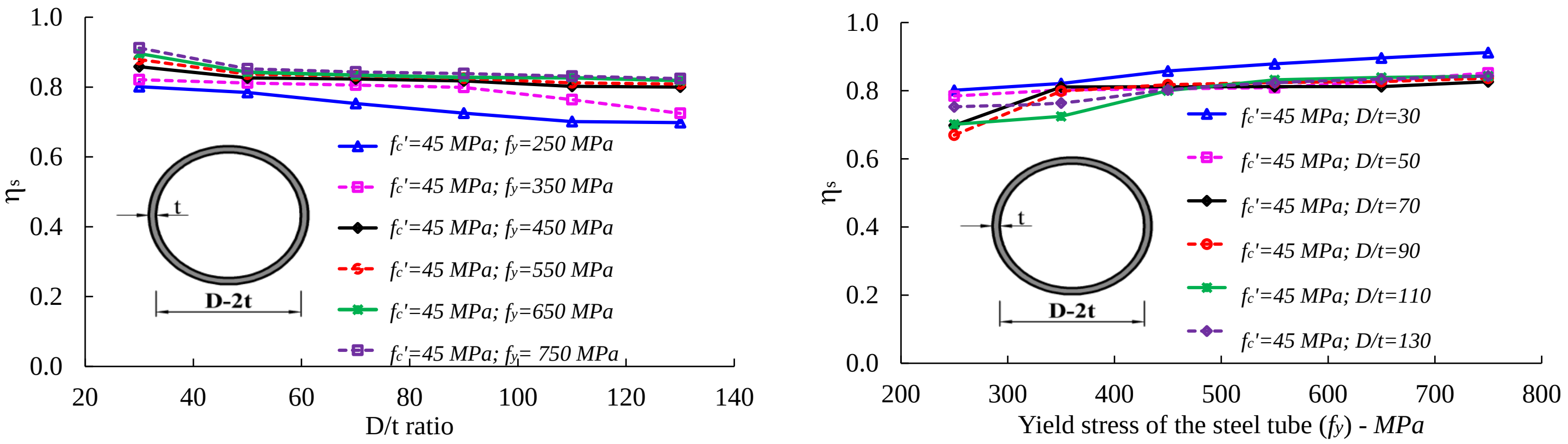

Fig. 15 Analyze the effects of parameters D/t ratio and $f_y$ on $\eta_s$

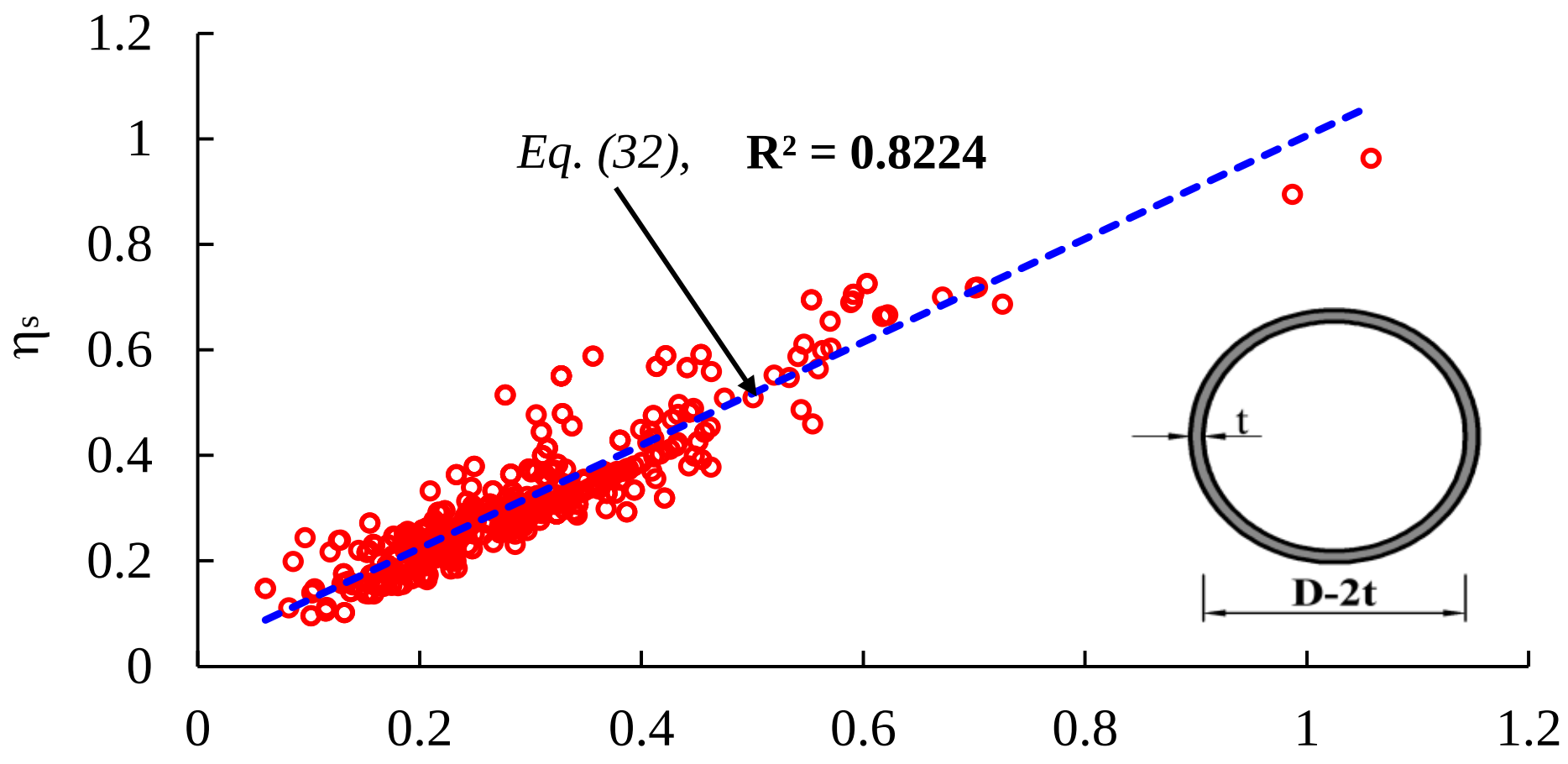

Fig. 16 Predict the impact coefficient $\eta_s$ of the steel tubes

The intensification factor $\eta_c$ for the concrete is mainly affected by, D/t or $f_c'$, as shown in Fig. 17. Increasing decreasing D/t or $f_c'$ leads to increasing concrete confinement. As can be seen from Fig. 17, $\eta_c$ decreases dramatically as D/t or $f_c'$ increases initially. Then $\eta_c$ becomes stable at a value larger than 1. In contrast, $\eta_c$ increases almost linearly with increasing D/t ratio, as shown in Fig. 18. Regression analysis indicates that $\eta_c$ may be expressed as a function of D/t ratio only. However, if *D/t* or $f_c'$ are introduced as additional terms, a better model can be produced for $\eta_c$ as shown in Fig. 18, where the value of $R^2$ is 0.9791. The equation to predict $\eta_c$ is given as follows:

$$\eta_c = 0.85 + 0.3 \times \left[\frac{D}{t}\right]^{0.328} \times \left(f_c'\right)^{0.1} \times \xi_c \tag{33}$$

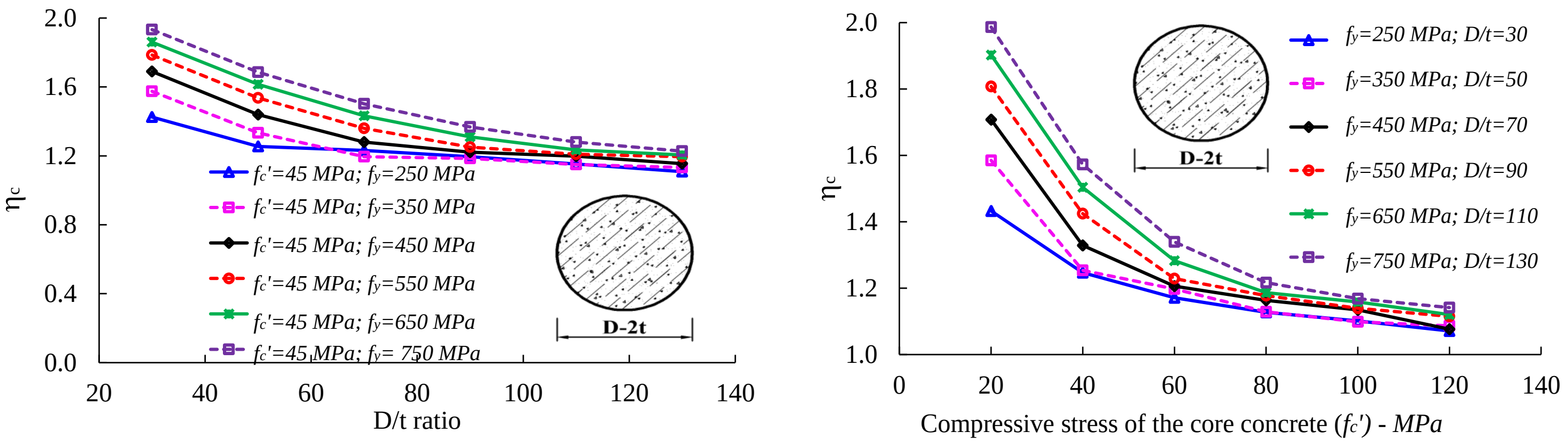

Fig. 17 Analyze the effects of parameters D/t ratio and $f_c{}'$ on $\eta_c$

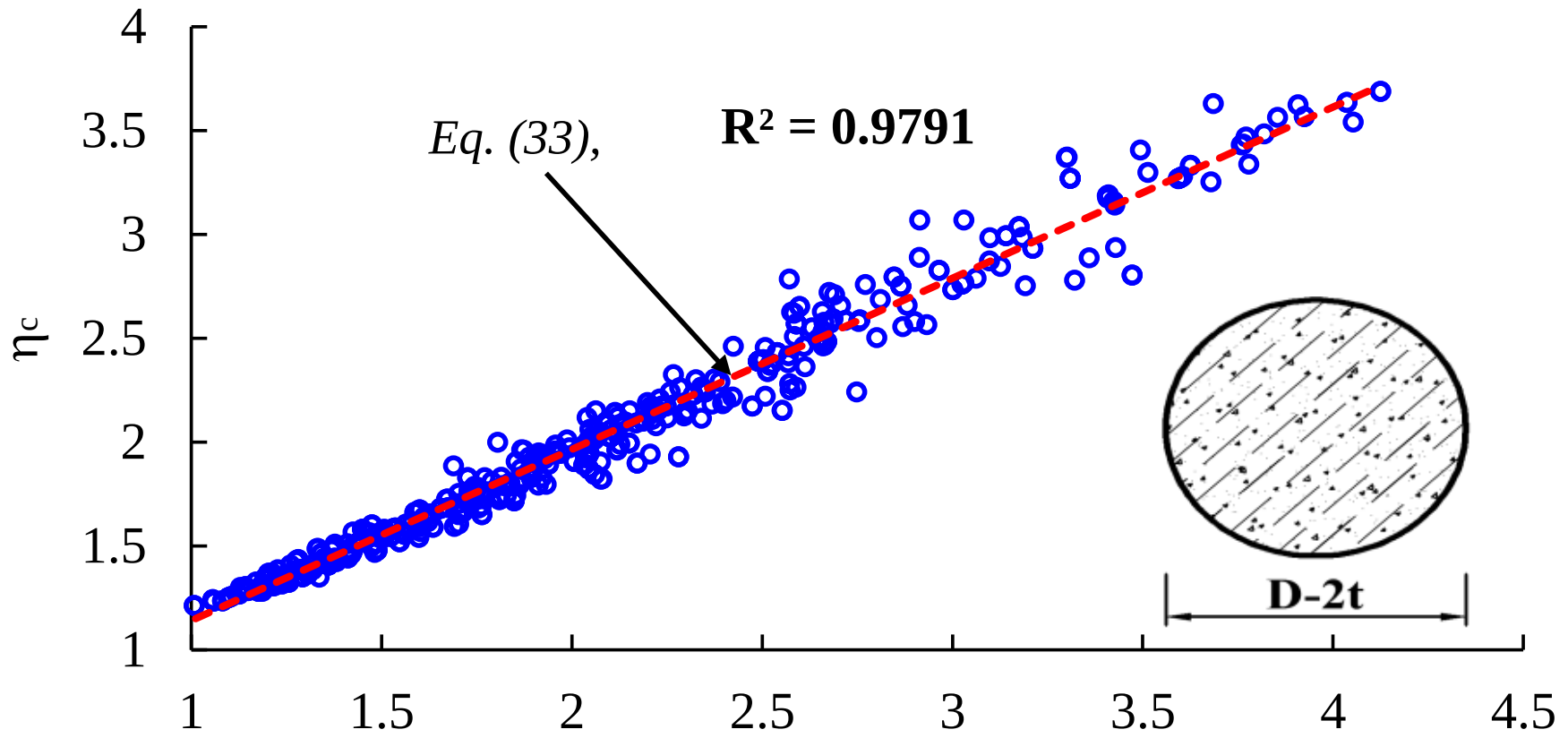

Fig. 18 Predict the impact coefficient $\eta_c$ of the core concrete

*5.4 Verification of the calculation formula*

The introduction and applying Equation (30) to 663 experimental, which the purpose of comparing the bearing capacity of the CFST columns with current design codes and previous formulas. The new design proposed formula improves the limits specified in design codes and evaluation of the previous formulas.

Table 7 Comparisons between measured and calculated results of CFST specimens.

| FE model | N | Mean | STD | CoV |
|---|---|---|---|---|
| **Current design codes** | | | | |
| Proposed design | 663/663 | 1.017 | 0.174 | 0.171 |
| EC4 2004 [18] | 472/663 | 1.13 | 0.17 | 0.15 |
| AISC 360-10 2010 [19] | 330/663 | 1.33 | 0.19 | 0.14 |
| CISC 2007 [21] | 663/663 | 1.18 | 0.39 | 0.33 |
| DBJ 13-51-2010 [22] | 663/663 | 1.06 | 0.21 | 0.20 |
| ACI 318 [20] | 279/663 | 1.03 | 0.12 | 0.12 |
| **Formulas design previous** | | | | |
| Proposed design | 663/663 | 1.017 | 0.174 | 0.171 |
| O'Shea and Bridge [1] | 580/663 | 0.77 | 0.22 | 0.28 |
| Yu et al [88] | 141/663 | 1.01 | 0.11 | 0.11 |
| Liu et al [89] | 619/663 | 0.98 | 0.12 | 0.13 |
| Sun [90] | 663/663 | 0.93 | 0.45 | 0.48 |
| Zhong and Miao [44] | 663/663 | 1.12 | 0.21 | 0.18 |
| Guo [91] | 553/663 | 0.81 | 0.10 | 0.12 |
| De Oliveira et al [92] | 629/663 | 1.27 | 0.21 | 0.17 |

Table 7, Fig. 19, Fig. 20 indicates that the Mean, STD, Cov are 1.017, 0.174, 0.171 respectively for the current design formulas. The calculated results by Eq (28) are slightly lower than the measure experimental. The predictions ultimate load design codes comparison with Eq (28) for normal strength ($f_c \leq 60$ MPa, $f_y \leq 460$ MPa) of columns, which there is still a difference of about ±15%. Besides, the current standards give more accurate predictions than previous prediction formulas. Moreover, when using the high and ultra-high-strength ($f_c$ >120 MPa, $f_y$ >550 MPa) in circular short columns, the prediction EC4 2004 [18], Yu et al [88], CISC 2007 [21], DBJ 13-51-2010 [22], Liu et al [89], Zhong and Miao [44] estimates give better results and allowable deviations about under ±10%. The prediction of the ultimate of AISC 360-10 2010 [19], ACI 318 [20], O'Shea and Bridge [1], Sun [90], Guo [91], De Oliveira et al [92] are large deviation when compared to the measured experiment, so the estimation of the predictions is not reliable. Setting limits in current standards is limited in predicting the bearing capacity of CFST columns when using high-strength materials. On the other hand, the proposed models previously were ignored effect geometry imperfections, residual stresses in tube steel, which leads to an assessment of the bearing capacity are limited.

Lastly, the results of the comparison between the measured and calculated axial strengths of all specimens when applying on Eq (30) are reasonable and accurate. The prediction ultimate load for short CFST columns through the contribution factors ηs and ηc from the simulation results are considered initial imperfection and residual stress, which leads to increased accuracy in the prediction.

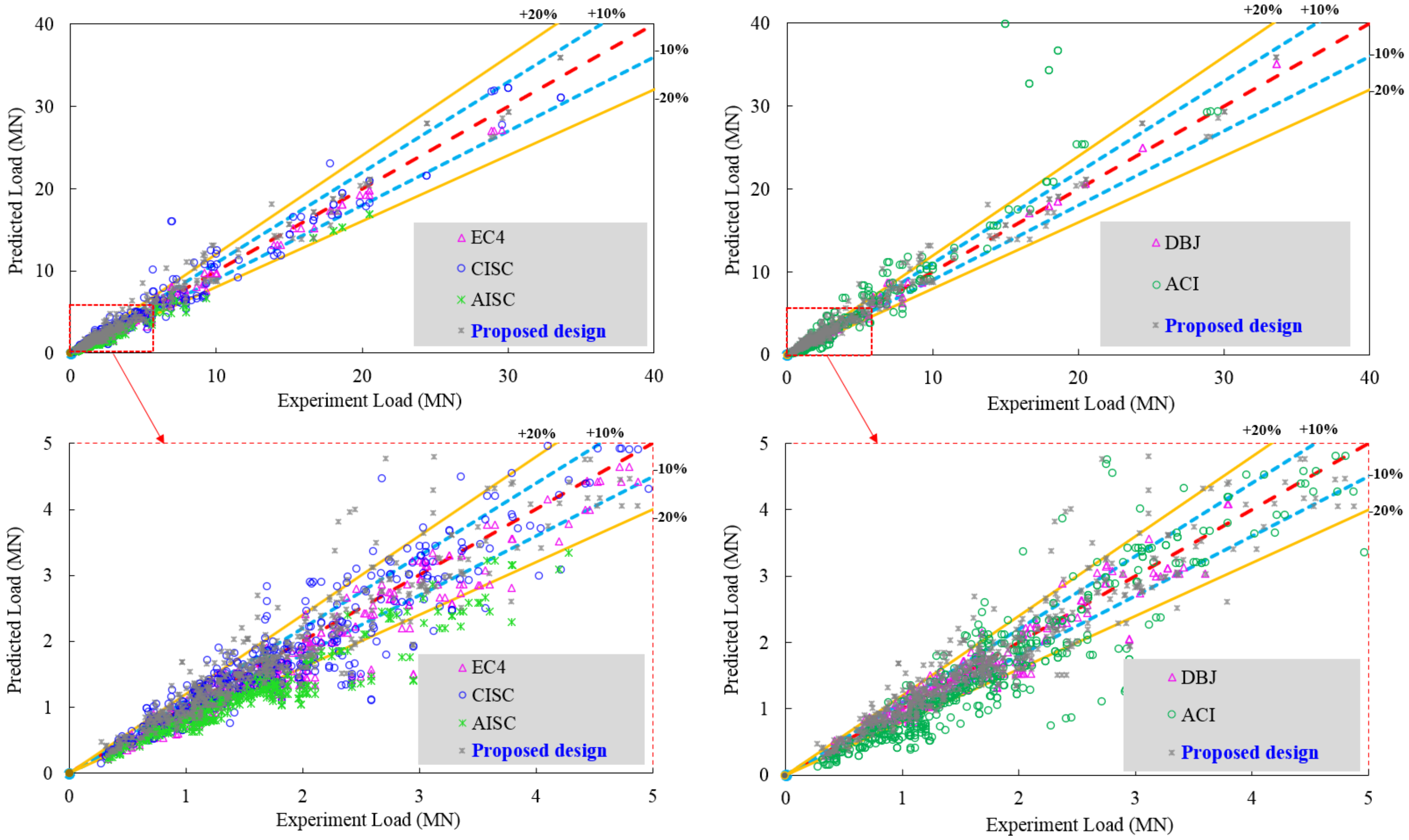

Fig. 19 Comparison between the proposed design and current design codes

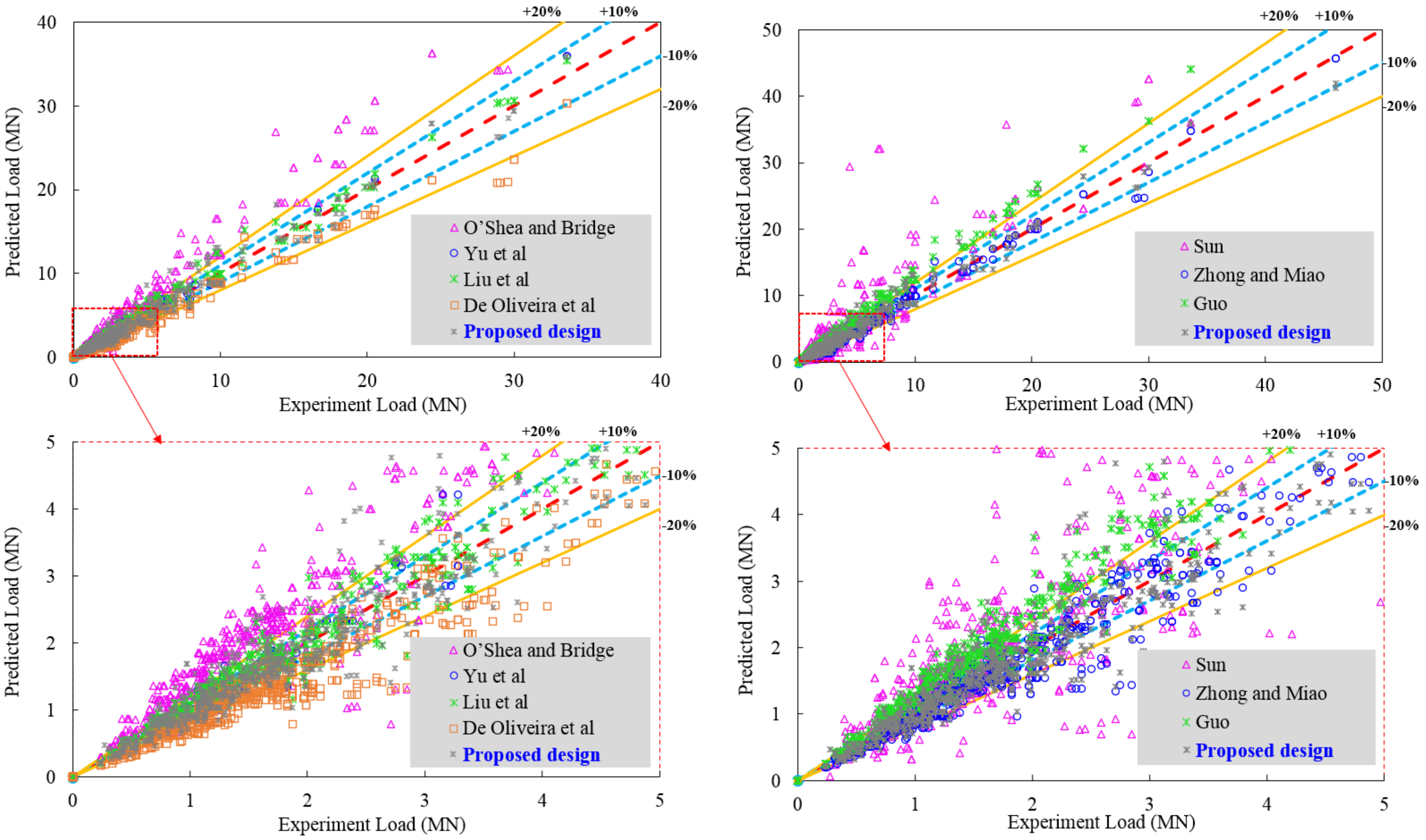

Fig. 20 Comparison between the proposed previous formula and proposed design

## 6.Conclusion

The authors proposed a new finite element model employing ABAQUS for predicting accurately the static behavior of concrete-filled steel circular tube column (CFSCT) subjected axial compression. Some conclusions can be drawn as follows:

- To optimize and minimize the computational time of the choosing axisymmetric model is more appropriate than the full 3D model in CFST columns simulation.
- Consider initial imperfections and residual stress in steel tubes short circular columns
- Determine the formula for calculating the confinement stress in a CFST column under axial compression by the regression method.

- Comparing with Han et al.’s and Hu et al.’s FE models, the proposed model predicts more exactly the ultimate strength of CFSCT columns.
- The proposed model can capture accurately the static behavior of CFSCT columns more than previous studies.
- The current standards give more accurate predictions than previous prediction formulas when using normal material. But current standard can't predictions for high and ultra high strength, because of the limited rules.
- The prediction ultimate load for short CFST columns through the contribution factors $\eta_s$ and $\eta_c$ from the simulation results are considered initial imperfection and residual stress, which leads to increased accuracy in the prediction.

## Acknowledgments

This work was supported by a grant from the Human Resources Development of the Korea Institute of Energy Technology Evaluation & Planning (KETEP) funded by the Korea government Ministry of Knowledge Economy (No. 20104010100520) and by the National Research Foundation of Korea (NRF) grant funded by the Korea government (MEST) (No. 2011-0030847).